\input ./style/arxiv-general.cfg
\documentclass[MSNbibl,number,citesort,seceqn,dvips]{arxbj}
\makeatletter
   \@ifpackageloaded{graphicx}{}{\usepackage{graphicx}}
\makeatother
\usepackage{multirow}

\volume{22}
\issue{2}
\pubyear{2016}
\firstpage{1093}
\lastpage{1112}
\doi{10.3150/14-BEJ687}% Updated by VTEXPTS2LaTeX.exe, 10.12.2014 10:24
\docsubty{FLA}

\makeatletter
\newcommand{\rrvert}{\vert}
\newcommand{\rrVert}{\Vert}
\newcommand{\llvert}{\vert}
\newcommand{\llVert}{\Vert}
\renewcommand{\mid}{|}
\def\R{{\mathbb R}}
\def\iny{\infty}

\newtheorem{lemma}{Lemma}[section]
\newtheorem{theorem}{Theorem}[section]
\newtheorem{corollary}{Corollary}[section]
\makeatother

\begin{document}
\begin{frontmatter}

\title{Behavior of R-estimators under measurement~errors}
\runtitle{R-estimates in ME models}

\begin{aug}
%%%% inicialai - be tarpu
% Corresponding author: Jana Jureckova - jurecko@karlin.mff.cuni.cz% Updated by VTEXPTS2LaTeX.exe, 15.12.2014 14:47
%Updated by VTEXPTS2LaTeX.exe, 10.12.2014 10:24
\author[A]{\inits{J.}\fnms{Jana}~\snm{Jure\v{c}kov\'a}\corref{}\thanksref{A,e1,u1}\ead[label=e1,mark]{jurecko@karlin.mff.cuni.cz}\ead[label=u1,url,mark]{http://www.karlin.mff.cuni.cz/\textasciitilde jurecko}},
\author[B]{\inits{H.L.}\fnms{Hira L.}~\snm{Koul}\thanksref{B}\ead[label=e3]{koul@stt.msu.edu}},
\author[D]{\inits{R.}\fnms{Radim}~\snm{Navr\'atil}\thanksref{A,D,e2}\ead[label=e2,mark]{navratil@karlin.mff.cuni.cz}}
\and
\author[C]{\inits{J.}\fnms{Jan}~\snm{Picek}\thanksref{C}\ead
[label=e4]{jan.picek@tul.cz}}
%%\runauthor{} %% autoJan
%\dedicated{}
\address[A]{Faculty of Mathematics and Physics, Charles University,
Sokolovsk\'{a} 83, CZ-186 75 Prague 8,\\ Czech Republic. \printead{e1,e2};\\ \printead{{u1}}}
\address[B]{Statistics and Probability, Michigan State University,
A435 Wells Hall, East Lansing, Michigan, USA.\\ \printead{e3}}
\address[D]{Department of Mathematics and Statistics, Masaryk
University Brno, Kotl\'{a}\v{r}sk\'{a} 2, 611 37 Brno,\\ Czech
Republic}
\address[C]{Applied Mathematics, Technical University, Voron\v{e}\v
{z}sk\'a 13, Liberec, Czech Republic.\\ \printead{e4}}
\end{aug}

% HISTORY:
%
\received{\smonth{2} \syear{2014}}% Updated by VTEXPTS2LaTeX.exe,
%10.12.2014 10:24
%
\revised{\smonth{10} \syear{2014}}% Updated by VTEXPTS2LaTeX.exe,
%10.12.2014 10:24

% ABSTRACT
%
\begin{abstract}
As was shown recently, the measurement errors in regressors affect only
the power of the rank test,
but not its critical region. Noting that, we study the effect of
measurement errors on R-estimators
in linear model. It is demonstrated that while an R-estimator admits a
local asymptotic bias,
its bias surprisingly depends only on the precision of measurements and
does neither depend
on the chosen rank test score-generating function nor on the
{regression model error distribution}.
The R-estimators are numerically illustrated and compared with the LSE
and $L_1$ estimators in this situation.
\end{abstract}

% KEYWORDS
% visi is mazosios raides ir pagal abecele
%
\begin{keyword}
\kwd{contiguity}
\kwd{linear rank statistic}
\kwd{linear regression model}
\kwd{local asymptotic bias}
\kwd{measurement error}
\kwd{R-estimate}
\end{keyword}
\end{frontmatter}

%s1 #&#
\section{Introduction}\label{sec1}
Measurement technologies are often affected by random errors; if the
goal of the experiment is to estimate a parameter, then the estimate is
biased, and thus inconsistent. This problem appears in the analytic
chemistry, in environmental monitoring, in modeling astronomical data,
in biometrics, and practically in all parts of the reality. Moreover,
some observations can be undetected, for example, when the measured
flux (light, magnetic) in the experiment falls below some flux limit.
In econometrics, the errors can be a result of misreporting by
subjects, miscoding by the collectors of the data, or by incorrect
transformation from initial reports. An essential part
of measuring techniques, used, for example, in the analytic chemistry,
is the construction of a calibration curve -- the result for an unknown
sample is then determined by interpolation.
Robust calibration methods were developed in
\cite{Ivo}. However, even the
calibration can be affected by measurement errors.
The mismeasurements make the statistical inference biased, and they
distort the trends in the data.

A variety of functional models have been proposed for handling
measurement errors
in regression models. Either the regressor or the response or both can
be affected by random errors. Technicians, geologists and other
specialists are aware of this problem, and try to reduce the bias with
various ad hoc procedures. The bias cannot be completely eliminated or
substantially reduced unless we have some additional knowledge on the
behavior of measurement errors.
The papers dealing with practical aspects of measurement error models
include \cite{Akritas,Hyk,Kelly,Marques,Rocke},
among others.

Adcock \cite{Adcock} was probably the first to realize the importance
of the situation. There exists a rich literature on the statistical
inference in the error-in-variables (EV) models as is evidenced by the
monographs of Fuller \cite{Fuller1987}, Carroll et al. \cite
{Carroll2006}, and Cheng and Van Ness \cite{vanNess1999}, and the
references therein. The monographs \cite{Fuller1987} and \cite
{vanNess1999} deal mostly with classical Gaussian set up while \cite
{Carroll2006} discusses numerous inference procedure under
semi-parametric set up. Nonparametric methods in EV models are
considered in \cite{Carroll1999,Carroll2007} and in
references therein, and in \cite{FanTruong1993}, among others. The
regression quantile theory in the area of EV models was started by He
and Liang \cite{HeLiang}. Arias, Hallock and Sosa-Escudero
\cite{Arias2001} used an instrumental
variable estimator for quantile regression, considering biases arising
from unmeasured ability and measurement errors. The problem of
mismeasurement is also of interest in the econometric literature: \cite
{Hausman2001} and \cite{Hyslop} described the recent developments in
treating the effect of mismeasurement on econometric models.

The advantage of \textit{rank and signed rank procedures} in the
measurement errors models was discovered recently in \cite{Jurecko1,NavratilSaleh,Navratil,Saleh2} and in \cite
{JISA}; the latter made
a detailed analysis of rank procedures in the linear model with a
nonlinear nuisance regressor and under various kinds of measurement
errors. Namely the rank tests can be recommended in this situation: it
is shown in \cite{Jurecko1} that the critical region of the rank test
for regression is insensitive to measurement errors in regressors under
very general conditions; the errors affect only the power of the test.
However, against expectations following from the invariance of the
ranks, due to which an estimate of a nuisance parameter in \cite
{Jurecko1} was consistent for every fixed value of the same, we show
that the R-estimator of slope parameter $\bolds\beta$ in linear model
is biased.
More precisely, we show that, unless $\bolds\beta=\mathbf0$, the
R-estimator is biased even in a local neighborhood of $\mathbf0$.
Hence,
we cannot have an unbiased estimator of any kind in this situation,
unless we have some additional information on the measurement errors.

As we further show in the present paper,
{surprisingly} the local asymptotic bias of R-estimators
neither depend{s}
on the chosen rank test score-generating function{s} nor on the
unknown distribution of the model errors. It depends only on value of
{slope parameter vector}
and on the covariance matrix of the measurement error {distribution} of
regressors.

%s2 #&#
\section{Model and preliminary considerations}\label{sec2}
Consider the linear regression model
%
%e2.1 #&#
%
\begin{equation}
\label{2} Y_{ni}=\beta_0+{\mathbf x}_{ni}^{\top}{
\bolds\beta}+e_{ni}, \qquad i=1,\ldots,n
\end{equation}
with unknown parameters $\beta_0\in\mathbb{R}^1$, $\bolds{\beta
}\in\mathbb{R}^{p}$.
The regressors ${\mathbf x}_{ni}$ are either deterministic or random
and affected by additive random measurement errors,
so that instead of ${\mathbf x}_{ni}$ we observe ${\mathbf
w}_{ni}={\mathbf x}_{ni}+{\mathbf v}_{ni}, i=1,\ldots,n$,
where ${\mathbf v}_{n1},\ldots,{\mathbf v}_{nn}$ are $p$-dimensional
random errors, identically distributed with an unknown distribution,
and independent of the errors $e_{ni}, 1\le i\le n$. Moreover, there
are additive measurement errors in the responses, thus instead of
$Y_{ni}$ we observe $Y_{ni}^*=Y_{ni}+u_{ni}$, where $u_{n1},\ldots
,u_{nn}$ are i.i.d. random variables. Thus in terms of the observable
responses and predicting variables, our regression model becomes
%
%e2.2 #&#
%
\begin{equation}
\label{contaminated} Y_{ni}^*=\beta_0+{\mathbf
w}_{ni}^{\top}{\bolds\beta}+e_{ni}^*, \qquad i=1,\ldots,n,
\end{equation}
where $e_{ni}^*=e_{ni}^*(\bolds{\beta})=e_{ni}+u_{ni}-\mathbf
{v}_{ni}^{\top}\bolds{\beta}, i=1,\ldots,n$ are i.i.d random variables.

We are interested in R-estimator of the slope vector $\bolds\beta$,
considering $\beta_0$ as nuisance parameter.
To define these estimators, let ${R}_{ni}(\mathbf b)$
be the rank of the residual
\begin{eqnarray*}
Y_{ni}^*-\mathbf w_{ni}^{\top}\mathbf
b&=&e_{ni}+u_{ni}+\mathbf x_{ni}^{\top}
\bolds\beta-\mathbf w_{ni}^{\top}\mathbf b
\\
&=&e_{ni}+u_{ni}-{\mathbf w}_{ni}^{\top}{
\mathbf b^*}-\mathbf v_{ni}^{\top}\bolds\beta, \qquad i=1,\ldots,n,
\end{eqnarray*}
where $\mathbf b^*=\mathbf b-\bolds\beta$.
We shall work with the vector of linear rank statistics
%
%e2.3 #&#
%
\begin{equation}
\label{S5} \mathbf S_n(\mathbf b)= { \bigl(S_{nj}(
\mathbf b); j=1,\ldots,p \bigr)^{\top}} =n^{-1/2}\sum
_{i=1}^n ({\mathbf w}_{ni}-\bar{\mathbf
w}_n)a_n\bigl({R}_{ni}(\mathbf b)\bigr),
\end{equation}
where the scores $a_n(i), 1\le i\le n$ are nondecreasing in $i$ and
$\sum_{i=1}^n a_n(i)=0$.

Hodges and Lehmann \cite{Hodges} introduced a class of estimators of
the location parameter $\theta$ in one- and two-sample location
models, by inverting a class of rank tests for $\theta$. This
methodology was extended to linear regression models without
measurement error by Jure\v{c}kov\'a \cite{Jur1971}, where an
estimator of $\bolds\beta$ is defined as
\[
\widehat{\bolds\beta}_n=\mathop{\arg\min}_{\mathbf b\in\R_p}\sum
_{j=1}^p\bigl\llvert S_{nj}(
\mathbf b)\bigr\rrvert. %
\]
This estimator can be seen to be asymptotically equivalent to an
estimator obtained by inverting the equations $S_{nj}(\mathbf b)=0,
j=1,\ldots, p$. Note that this latter estimator is precisely an
extension of the Hodges--Lehmann estimator from one- and two-sample
location models to linear regression models without measurement error.
Under more general conditions, the
R-estimators are studied by Koul \cite{Koul2000}.

On the other hand, Jaeckel \cite{Jaeckel} called an analog of the function
%
%e2.4 #&#
%
\begin{eqnarray}
\label{S6} \mathcal D_n(\mathbf{b})&=& \sum
_{i=1}^n \bigl[Y_{ni}^*-
\mathbf{w}_{ni}^{\top}\mathbf{b} \bigr]\bigl(a_n
\bigl({R}_{ni}(\mathbf{b})\bigr)-\bar{a}_n\bigr),
\end{eqnarray}
as a measure of rank dispersion of residuals, in the case of no
measurement error where $\mathbf{w}_{ni}$'s are replaced by $\mathbf
{x}_{ni}$'s. He\vspace*{1pt} showed that $\mathcal D_n(\mathbf{b})$ is convex and
piecewise linear in $\mathbf b\in\mathbb R^p$. He also showed that
$-n^{1/2}\mathbf{S}_n(\mathbf{b})$ is the subgradient of $\mathcal
D_n(\mathbf{b})$; hence the estimator defined as a minimizer of
$\mathcal D_n$ exists and is equivalent to the above estimators based
on $\mathbf S_n$.
Both of these estimators are asymptotically equivalent, and Jaeckel's
definition of R-estimator is now generally used in the literature. We
are using this definition of R-estimator throughout this paper.

In the absence of measurement errors, that is, if $\mathbf
w_{ni}=\mathbf x_{ni}, u_{ni}=0, i=1,\ldots,n$, the estimator
$\widehat{\bolds\beta}_n$ is consistent and asymptotically normal.
However, $\widehat{\bolds\beta}_n$ is biased in the presence of measurement
errors, even asymptotically, unless the true $\bolds\beta=\mathbf0$.
Furthermore, we show that it is even asymptotically locally biased in
the sense that the asymptotic distribution of $n^{1/2}(\widehat{\bolds
\beta}_n - n^{-1/2}\bolds\beta^0)$, with a fixed
$\bolds\beta^0\in\mathbb R^p$, converges to a normal distribution
with nonzero mean vector and some positive definite covariance matrix.

In the sequel, all limits are taken as $n\to\infty$, unless mentioned
otherwise, $\stackrel{p}{\rightarrow}$ denotes the convergence in
probability. We
shall now describe the needed assumptions on the underlying entities.
\begin{longlist}[(A.1)]
\item[(A.1)] The score generating function $\varphi\dvtx (0,1)\mapsto
\mathbb R$ is nondecreasing, square-integrable and skew-symmetric on
$(0,1)$, that is,
satisfies $\varphi(1-t)=-\varphi(t), 0<t<1$. The scores $a_n(i),
i=1,\ldots,n$ are generated by $\varphi$ in either of the following
two ways:
\[
a_n(i)=\varphi\biggl(\frac{i}{n+1} \biggr)\quad\mbox{or}\quad
a_n(i)=\mathbb E\varphi(U_{n:i} ),\qquad i=1,\ldots,n,
\]
where $U_{n:1}\leq\cdots\leq U_{n;n}$ are order statistics pertaining
to the sample of size $n$ from the uniform $(0,1)$ distribution.
\item[(F.1)] Distribution function $F$ of the model errors $e_{ni}$
has an absolutely continuous density $f$ with a.e. derivative $f'$. %,

\item[(F.2)] For every $u\in\R$, $\int
(\llvert f'(x-tu)\rrvert ^j/f^{j-1}(x)) \,\mathrm{d}x \to\int(\llvert
f'(x)\rrvert ^j/f^{j-1}(x))
\,\mathrm{d}x<\infty$, as $t\to0$, $j=2,3$.

\item[(V.1)] The measurement errors $\{u_{ni}, 1\le i\le n\}$ are
independent of $\{e_{ni}, \mathbf v_{ni}, 1\le i\le n\}$ and i.i.d.
with generally an unknown absolutely continuous density $h$, having
finite Fisher's information for location.
\item[(V.2)] The measurement error $\mathbf v_{ni}$ {is} independent
of $e_{ni}$ and {its} $p$-dimensional distribution function $G$ has a
continuous density $g$, generally unknown, $i=1,\ldots,n$.
\item[(V.3)]
$\mathbb E \mathbf V_n \rightarrow\mathbf V$ where\vspace*{1pt} $\mathbf
V_n=n^{-1}\sum_{i=1}^n (\mathbf v_{ni}-\bar{\mathbf v}_n)(\mathbf
v_{ni}-\bar{\mathbf v}_n)^{\top}$ and $\mathbf V$ is a positive
definite $p\times p$ matrix.
Moreover,
$\sup_{n\ge1}\mathbb E (\llVert \mathbf v_{n1}\rrVert ^3+ \llVert
\mathbf x_{n1}\rrVert
^3 ) <\infty$.
\item[(V.4)] $\mathbb E [n^{-1}\sum_{i=1}^n(\mathbf v_{ni}-\bar
{\mathbf v}_n)(\mathbf x_{ni}-
\bar{\mathbf x}_n)^{\top} ]\rightarrow\mathbf0$.
\item[(X.1)] If the regressors $\mathbf x_{ni}$ are nonrandom, then
assume that $\mathbf Q_n\to\mathbf Q$, where
\[
\mathbf Q_n=n^{-1} \sum_{i=1}^n
(\mathbf x_{ni}-\bar{\mathbf x}_n) (\mathbf
x_{ni}-\bar{\mathbf x}_n)^{\top},
\]
and $\mathbf Q$ is positive
definite $p\times p$ matrix. Moreover,
\[
\frac{1}n\max_{1\leq i\leq n}(\mathbf x_{ni}-\bar{
\mathbf x}_n)^{\top
}\mathbf(\mathbf Q_n)^{-1}(
\mathbf x_{ni}-\bar{\mathbf x}_n)\to0.
\]
\item[(X.2)] If the regressors $\mathbf x_{ni}$ are random, then
assume that they are independent of $e_{ni}, u_{ni}, \mathbf
v_{ni}$, $i=1,\ldots,n$, and
\[
\mathbb E \Biggl[n^{-1} \sum_{i=1}^n
(\mathbf x_{ni}-\bar{\mathbf x}_n) (\mathbf
x_{ni}-\bar{\mathbf x}_n)^{\top} \Biggr]\rightarrow
\mathbf Q,
\]
where $\mathbf Q$ is positive
definite $p\times p$ matrix.
\end{longlist}

Let $m(\cdot), M(\cdot)$ denote the density and distribution
function of $e_{ni}+u_{ni}, i=1,\ldots,n$, that is,
$m(z)=\int f(z-t)h(t)\,\mathrm{d}t$. The density is absolutely continuous and has
finite Fisher's information $\mathcal I(m)$.
We need to define
%
%e2.5 #&#
%
\begin{eqnarray}
\label{m}  \gamma_m &=&-\int_{\mathbb R^1}\varphi
\bigl(M(z)\bigr)\,\mathrm{d}m(z), \qquad A_m^2(\varphi)=
\gamma_m^{-2}\int_0^1
\varphi^2(u) \,\mathrm{d}u,
\nonumber\\[-8pt]\\[-8pt]\nonumber
\mathbf B&=&-(\mathbf Q+\mathbf V)^{-1}\mathbf V\bolds
\beta^0.
\nonumber
\end{eqnarray}
The following theorem gives the asymptotic distribution of the
estimator $\widehat{\bolds\beta}_n$ when the
true parameter value is
%
%e2.6 #&#
%
\begin{equation}
\label{alt} \bolds\beta_n=n^{-1/2}\bolds
\beta^0,\qquad \bolds\beta^0\in\mathbb R^p\mbox{ fixed}.
\end{equation}

\begin{theorem}\label{Theorem0}
{Assume} the conditions \textup{(A.1)}, \textup{(F.1)}--\textup{(F.2)},
\textup{(V.1)}--\textup{(V.4)}, \textup{(X.1)}--\textup{(X.2)} {hold}.
{When the true parameter value\vspace*{1pt} is $\bolds\beta_n$,}
the R-estimator $\widehat{\bolds\beta}_n$ is asymptotically normally
distributed with the bias $\mathbf B=-(\mathbf Q+\mathbf
V)^{-1}\mathbf V\bolds\beta^0$, that is,
%
%e2.7 #&#
%
\begin{eqnarray}
\label{normality} n^{1/2}(\widehat{\bolds\beta}_n-\bolds
\beta_n)\stackrel{\mathcal D} {\rightarrow}\mathcal N_p
\bigl(\mathbf B,(\mathbf Q+\mathbf V)^{-1} A_m^2(
\varphi) \bigr).
\end{eqnarray}
\end{theorem}

Theorem \ref{Theorem0} will be proved in several steps; the proof is
given in Section~\ref{sec3}. The numerical illustrations of the results are
given in subsequent Section~\ref{sec4}.

\begin{corollary}\label{normality2}
Under the conditions of Theorem \ref{Theorem0} {and under $\bolds
\beta=\bolds\beta_n=n^{-1/2}\bolds\beta^0$,} the \mbox{R-}estimator
$\widehat{\bolds\beta}_n$
has asymptotic normal distribution
%
%e2.8 #&#
%
\begin{equation}
\label{corollary} n^{1/2}\bigl(\widehat{\bolds\beta}_n-(\mathbf
Q+\mathbf V)^{-1}\mathbf Q \bolds\beta_n\bigr)\stackrel{
\mathcal D} {\rightarrow}\mathcal N_p \bigl(\mathbf0,(\mathbf Q+
\mathbf V)^{-1} A_m^2(\varphi) \bigr).
\end{equation}
\end{corollary}

Notice that the local asymptotic bias cannot be controlled by the
choice of the score-generating function $\varphi$;
this choice can only influence the asymptotic variance factor of the
estimator. The magnitude of the bias fully depends
on the precision of the measurements, namely on the matrix~$\mathbf V$.
The measurement errors in the responses $Y_{ni}$
affect only the asymptotic variance, not the bias. The result is
entirely nonparametric, valid for
classes of distributions of model {and measurement} errors, demanding
only finite first moment and finite (and
positive) Fisher's information for location of the model error
distributions, and finite third moment for
measurement error distributions.

Consider the two measurement methods with the same regressors (random
or nonrandom), with the respective limiting
covariance matrices $\mathbf V_1, \mathbf V_2$. Comparing the biases
in (\ref{normality}) for $\mathbf V_1$ and
$\mathbf V_2$, the first method is considered being more precise than
the second one if the matrix
$(\mathbf V_2+\mathbf Q)^{-1}\prec(\mathbf V_1+\mathbf Q)^{-1}$;
otherwise speaking, if $\mathbf Q^{-1}\mathbf
V_1\prec\mathbf Q^{-1}\mathbf V_2${,}
where {the} ordering $\mathbf A\prec\mathbf B$ means that $\mathbf
B-\mathbf A$ is {a} positive definite
{matrix}.

%s3 #&#
\section{Proof of Theorem \texorpdfstring{\protect\ref{Theorem0}}{2.1}}\label{sec3}
We shall prove Theorem \ref{Theorem0} in several steps. Notice that if
we observe $Y_{ni}^*=Y_{ni}+u_{ni}$ instead of
$Y_{ni}$, then $e_{ni}^*=e_{ni}+u_{ni}, i=1,\ldots,n$ are still
i.i.d. random variables {with density $m(z)=\int
f(z-t)h(t)\,\mathrm{d}t$. The steps of the proof are parallel for both densities
$f$ and $m$ of model errors; measurement errors in
the $Y_{ni}$ affect only the asymptotic variance of the estimate, not
the bias. Noting this, we shall prove the theorem
assuming $u_{ni}\equiv0. i=1,\ldots,n$.}
In the sequel, we shall suppress the subscript $n$ whenever it does not
cause a confusion.

The steps of the proof are as follows:
\begin{longlist}[(2)]
\item[(1)] Asymptotic representation of the linear rank statistic
%
%e3.1 #&#
%
\begin{equation}
\label{S5a} \mathbf S_n(\mathbf0,\mathbf0)=n^{-1/2}\sum
_{i=1}^n ({\mathbf w}_{ni}-\bar{
\mathbf w}_n)a_n\bigl({R}_{ni}(\mathbf0)\bigr)
\end{equation}
with the sum of independent summands. {Here $\mathbf w_{ni}=\mathbf
x_{ni}+\mathbf v_{ni}, i=1,\ldots,n$, while
$\mathbf x_{n1},\ldots, \mathbf x_{nn}$ are either i.i.d. random
vectors or nonrandom vectors, and $\mathbf
v_{n1},\ldots,\mathbf v_{nn}$ are i.i.d. random vectors.}
\item[(2)] Contiguity\vspace*{1pt} of the sequence $\{ Q_n\}$ of distributions of
$(e_{ni}-(\mathbf w_{ni}-\bar{\mathbf
w}_n)^{\top}\mathbf b_n^*-(\mathbf v_{ni}-\bar{\mathbf v}_n)^{\top
}\bolds\beta_n)$, with $\mathbf
b_n^*=n^{-1/2}\mathbf b^0, \bolds\beta_n=n^{-1/2}\bolds\beta^0$ for
$\mathbf b^0, \bolds\beta^0 \in\mathbb R^p$ fixed, with respect
to the sequence $\{P_n\}$ of distributions of
$e_{ni}$, $i=1,\ldots,n$.
\item[(3)] Asymptotic representation of the linear rank statistic
(\ref{S5}) under contiguous sequence of
distribution $\{ Q_n\}$, and the resulting asymptotic linearity of
(\ref{S5}) in parameters $\mathbf b^0,
\bolds\beta^0$.
\item[(4)] Uniform asymptotic quadraticity of $\mathcal D_n$ in
parameters $\mathbf b^0, \bolds\beta^0$ under
$\{ Q_n\}$, as a result of (3) and of the convexity of
$\mathcal D_n$.
\item[(5)] Resulting asymptotic distribution and bias of $\widehat
{\bolds\beta}_n$ in the case $u_{ni}\equiv0,
i=1,\ldots,n$.
\item[(6)] Asymptotic distribution and bias of $\widehat{\bolds\beta
}_n$ in the case of nonzero $u_{ni},
i=1,\ldots,n$.
\end{longlist}

%s3.1 #&#
\subsection{Asymptotic representation of $\mathbf{S}_n(\mathbf{0},\mathbf{0})$}\label{sec3.1}
Assume
that $u_{ni}=0, i=1,\ldots,n$. That is, for now we assume that
there is no measurement
error in the response variables $Y_{ni}$.
Let
\begin{eqnarray*}
\mathbf Z_n=n^{-1/2}\sum_{i=1}^n(
\mathbf w_{ni}-\bar{\mathbf w}_n)\varphi
\bigl(F(e_{ni})\bigr).
\end{eqnarray*}
We are ready to state and prove
the following lemma.

\begin{lemma}\label{Lemma1} Under the conditions of Theorem \ref
{Theorem0}, the statistic $\mathbf S_n(\mathbf0, \mathbf0)$ admits
the asymptotic representation
%
%e3.2 #&#
%
\begin{equation}
\label{P52} \mathbf S_n(\mathbf0,\mathbf0)=\mathbf
Z_n+\mathrm{o}_p(1).
\end{equation}
\end{lemma}

\begin{pf}
The proof is adapted from \cite{Picek}. If $\mathbf
b=\bolds\beta=\mathbf0$, then
$(Y_{n1},\ldots,Y_{nn})=(e_{n1},\ldots,e_{nn})$. Let $R_{n1},\ldots
,R_{nn}$ denote their ranks. Further denote
$r_{ni}=a_n(R_{ni})-\varphi(F(e_{ni})), i=1,\ldots,n$.

Let $\sigma_j^2$ be the variance of {$w_{ij}, i=1,\ldots,n$,}
for $j=1,\ldots,p$, and let $s^2=\sum_{j=1}^p\sigma_j^2$.

Notice that $(r_{n1},\ldots,r_{nn})$ and {$(\mathbf w_1,\ldots
,\mathbf w_n)$} are independent. Consider the conditional
squared distance
\begin{eqnarray*}
&& \mathbb E_G \bigl\{(\mathbf S_n-\mathbf
Z_n)^{\top}(\mathbf S_n-\mathbf Z_n)
\mid e_1,\ldots,e_n \bigr\}
\\
&&\quad =n^{-1}\mathbb E_G \Biggl\{\sum
_{i=1}^n\sum_{k=1}^n(
\mathbf w_i-\bar{\mathbf w}_n)^{\top}(\mathbf
w_k-\bar{\mathbf w}_n)r_ir_k
\Big| e_1,\ldots,e_n \Biggr\}
\\
&&\quad = n^{-1}\sum_{i=1}^n\sum
_{k=1}^nr_ir_k
\mathbb E_G \Biggl\{\sum_{j=1}^p(w_{ij}-
\bar{w}_j) (w_{kj}-\bar{w}_j)\Big|
e_1,\ldots,e_n \Biggr\}
\\
&&\quad =n^{-1} \Biggl\{\sum_{i=1}^n
\sum_{k=1}^nr_ir_k
\sum_{j=1}^p(x_{ij}-
\bar{x}_j) (x_{kj}-\bar{x}_j)+s^2
\sum_{i=1}^n(r_i-
\bar{r})^2 \Biggr\}
\\
&&\quad =\sum_{j=1}^p \Biggl[n^{-1/2}
\sum_{i=1}^n(x_{ij}-\bar
{x}_j)r_i \Biggr]^2+s^2\sum
_{i=1}^n(r_i-
\bar{r})^2.
\end{eqnarray*}
Then (\ref{P52}) follows from \cite{Hajek1967} [Theorems V.1.4.a,b,
V.1.6.a].
\end{pf}%\hfill$\blacksquare$

%s3.2 #&#
\subsection{Contiguity}\label{sec3.2}
For any two probability measures $P$ and $Q$, absolutely continuous
with respect to a $\sigma$-finite measure $\nu$ with
$p=\mathrm{d}P/\mathrm{d}\nu$, $q=\mathrm{d}Q/\mathrm{d}\nu$, let
\begin{eqnarray*}
H(P,Q)= \biggl[\int(\sqrt{p}-\sqrt{q} )^2\,\mathrm{d}\mu\biggr]^{1/2}
= \biggl[2\int(1-\sqrt{pq} )\,\mathrm{d}\mu\biggr]^{1/2}
\end{eqnarray*}
denote the Hellinger distance between $P$ and $Q$.

Let $\{P_{n1},\ldots,P_{nn}\}$ and $\{Q_{n1},\ldots,Q_{nn}\}$ be two
triangular arrays of probability measures
defined on measurable space $({\mathcal X},{\mathcal A})$ with
densities $p_{ni}, q_{ni}$ with respect to
$\sigma$-finite
measures~$\mu_i$ [which can be also
$\mu_i=P_{ni}+Q_{ni}, i=1,\ldots,n$]. Denote $P_n^{(n)}=\prod_{i=1}^n
P_{ni}$ and
$Q_n^{(n)}=\prod_{i=1}^n Q_{ni}$ the product measures, $n=1,2,\ldots.$

Oosterhoff and van Zwet \cite{Oosterhoff79}
proved that $\{Q_n^{(n)}\}$ is contiguous with respect to $\{P_n^{(n)}\}$
if and only if
%
%e3.3 #&#
%e3.4 #&#
%
\begin{eqnarray}
\label{Zwet1} \limsup_{n\rightarrow\infty}\sum
_{i=1}^nH^2(P_{ni},Q_{ni})&<&
\infty,
\\
\label{Zwet2}  \lim_{n\rightarrow\infty}\sum_{i=1}^n
Q_{ni} \biggl\{\frac
{q_{ni}(X_{ni})}{p_{ni}(X_{ni})}\geq c_n \biggr\}&=&0 \qquad
\forall c_n\rightarrow\infty.
\end{eqnarray}
Note that in the case $P_{ni}\equiv P_n, p_{ni}\equiv p_n$, and
$Q_{ni}\equiv Q_n, q_{ni}\equiv q_n$, not depending on
$i$,
%
%e3.5 #&#
%
\begin{eqnarray}
\label{Zwet4} \sum_{i=1}^n
H^2(P_{ni},Q_{ni})&=&n\int\bigl[
\sqrt{q_n(z)}-\sqrt{p_n(z)} \bigr]^2\,\mathrm{d}z\nonumber
\\
&=& n\int
\frac{(q_n(z)-p_n(z))^2}{[\sqrt
{q_n(z)}+\sqrt{p_n(z)}]^2}\,\mathrm{d}z
\\
&\leq& n\int\frac{ (q_n(z)-p_n(z) )^2}{p_n(z)} \,\mathrm{d}z.\nonumber
\end{eqnarray}
Moreover, for $c_n>1$ and with $d_n=c_n-1$,
%
%e3.6 #&#
%
\begin{eqnarray}
\label{cont} \sum_{i=1}^n
Q_{ni} \biggl\{\frac{q_{ni}(X_{ni})}{p_{ni}(X_{ni})}\geq c_n \biggr\}&=&
nQ_n \biggl\{\frac{q_n(X_{n1})-p_n(X_{n1})}{p_n(X_{n1})}\geq d_n \biggr
\}\nonumber
\\
&\le& d_n^{-2}n \int\frac{\llvert q_n (x)-p_n(x)\rrvert ^2}{p_n^2(x)} q_n(x)
\,\mathrm{d}x
\nonumber\\[-8pt]\\[-8pt]\nonumber
&\le& d_n^{-2}n \int\frac{\llvert q_n (x)-p_n(x)\rrvert ^3}{p_n^2(x)}
\,\mathrm{d}x
\\
&&{} + d_n^{-2}n \int\frac{\llvert q_n
(x)-p_n(x)\rrvert ^2}{p_n(x)} \,\mathrm{d}x.
\nonumber
\end{eqnarray}

Now, let $Y_{ni}=\mathbf x_{ni}^{\top}\bolds\beta+e_{ni},
i=1,\ldots,n$. where the $e_{ni}$ are i.i.d. with
distribution function $F$ and density $f$, satisfying \textup{(F.1)} and
\textup{(F.2)}.
Consider the residuals
\begin{eqnarray}
\label{resid} Y_{ni}-(\mathbf w_{ni}-\bar{\mathbf
w}_n)^{\top}\mathbf b_n &=&e_{ni}+(
\mathbf x_{ni}-\bar{x}_n)^{\top}\bolds
\beta_n-(\mathbf w_{ni}-\bar{\mathbf w}_n)^{\top}
\mathbf b_n
\nonumber
\\
&=&e_{ni}-(\mathbf w_{ni}-\bar{\mathbf w}_n)^{\top}
\mathbf b_n^*-(\mathbf v_{ni}-\bar{\mathbf
v}_n)^{\top}\bolds\beta_n,
\nonumber
\end{eqnarray}
$i=1,\ldots,n$, where $\mathbf b_n=n^{-1/2}\mathbf b^0, \bolds
\beta_n=n^{-1/2}\bolds\beta^0, \mathbf
b_n^*=n^{-1/2}\mathbf b^{0*}, \mathbf b_0^*=\mathbf b^0-\bolds\beta
^0$, with fixed $\mathbf b^0,
\bolds\beta^0\in\mathbb R^p$.
Using (\ref{Zwet1}) and (\ref{Zwet2}), we shall prove the following lemma.

\begin{lemma}\label{Lemma2} Under the conditions
of Theorem \ref{Theorem0},
the sequence $\{Q_n^{(n)}\}$ is contiguous with respect to $\{
P_n^{(n)}\}$, where $Q_n^{(n)}=\prod_{i=1}^n Q_{ni}$, $P_n^{(n)}=\prod
_{i=1}^n P_{ni}$, where $P_{ni}$ is the distribution of $e_{ni}$ and
$Q_{ni}$ is the distribution of
$(e_{ni}-(\mathbf w_{ni}-\bar{\mathbf w}_n)^{\top}\mathbf
b_n^*-(\mathbf v_{ni}-\bar{\mathbf v}_n)^{\top}\bolds\beta_n),
i=1,\ldots,n$.
\end{lemma}

\begin{pf}
We shall distinguish the two cases: the $\mathbf x_{ni}$ are either
i.i.d. random vectors or nonrandom vector
components.

We start with the first case, where $\mathbf w_{n1},\ldots,\mathbf
w_{nn}$ are i.i.d. random vectors. Note that
$U_i:=(\mathbf w_{ni}-\bar{\mathbf w}_n)^{\top}\mathbf
b^{0*}+(\mathbf v_{ni}-\bar{\mathbf v}_n)^{\top}\bolds\beta^0,
i=1,\ldots,n$, are i.i.d. r.v.'s. Let {$k_1$} denote the common
density function of $U_i$. Then, $Q_{ni}, P_{ni}$ do
not depend on $i$ and
$q_n(x)\equiv\int f(x-n^{-1/2}u ) {k_1}(u) \,\mathrm{d}u, p_n(x)\equiv f(x)$.
Hence, by the Cauchy--Schwarz inequality, and the Fubini
theorem,
\begin{eqnarray*}
n\int\frac{ (q_n(x)-p_n(x) )^2}{p_n(x)} \,\mathrm{d}x&=& n\int\biggl\{\int\bigl
[f\bigl(x-n^{-1/2}u
\bigr)-f(x) \bigr]{k_1}(u)\,\mathrm{d}u \biggr\}^2\frac{\mathrm{d}x}{f(x)}
\\
&\le& n\int\!\!\int\bigl[f\bigl(x-n^{-1/2}u\bigr)- f(x)
\bigr]^2\frac
{{k_1}(u)}{f(x)}\,\mathrm{d}u \,\mathrm{d}x
\nonumber
\\
&\leq& n\int\!\!\int\biggl[\int_{-n^{-1/2}}^{n^{-1/2}} \bigl
\llvert u f^{\prime}(x-tu)\bigr\rrvert \,\mathrm{d}t \biggr]^2
\frac{{k_1}(u)}{f(x)}\,\mathrm{d}u \,\mathrm{d}x
\nonumber
\\
& \le& 2n^{1/2}\int\!\!\int\!\!\int_{-n^{-1/2}}^{n^{-1/2}}
\bigl\llvert f^{\prime}(x-tu)\bigr\rrvert^2\,\mathrm{d}t\, u^2
\frac{{k_1}(u)}{f(x)}\,\mathrm{d}u \,\mathrm{d}x
\nonumber
\\
& \le& 2n^{1/2}\int\!\!\int_{-n^{-1/2}}^{n^{-1/2}} \int
\frac{\llvert f^{\prime}(x-tu)\rrvert ^2}{f(x)} \,\mathrm{d}x\,u^2 {k_1}(u) \,\mathrm{d}u
\,\mathrm{d}t \qquad\forall n
\ge1.
\nonumber
\end{eqnarray*}
Hence, by (\ref{Zwet4}), \textup{(F.2)} applied with $j=2$, and by \textup{(V.3)}, which guaranteed $\int u^2{k_1}(u)\,\mathrm{d}u<\iny$,
%
%e3.7 #&#
%
\begin{eqnarray}
\label{zwet} \limsup_n \sum_{i=1}^n
H^2(P_{ni},Q_{ni})\leq2 I(f) \int
u^2{k_1}(u)\,\mathrm{d}u<\infty.
\end{eqnarray}

Similarly, the bound
\begin{eqnarray*}
n\int\frac{ (q_n(x)-p_n(x) )^3}{p_n^2(x)} \,\mathrm{d}x &\le& 2n^{1/2}\int\!\!\int
_{-n^{-1/2}}^{n^{-1/2}}
\int\frac{\llvert f^{\prime}(x-tu)\rrvert ^3}{f^2(x)} \,\mathrm{d}x \llvert
u\rrvert^3 {k_1}(u)
\,\mathrm{d}u \,\mathrm{d}t \qquad\forall n\ge1
\end{eqnarray*}
together with (\ref{cont}), \textup{(F.2)} applied with $j=3$, and
\textup{(V.3)}, which guaranteed $\int\llvert u\rrvert ^3 {k_1}(u)\,\mathrm{d}u<\iny$,
yield
\begin{eqnarray*}
&& \lim_n \sum_{i=1}^n
Q_{ni} \biggl\{\frac
{q_{ni}({Y}_{ni})}{p_{ni}({Y}_{ni})}\geq c_n \biggr\}
\\
&&\quad \le 2
\lim_n d_n^{-2} \biggl\{ \int\biggl(
\frac{\llvert f'(x)\rrvert }{f(x)} \biggr)^3 f(x)\,\mathrm{d}x \int\llvert
u\rrvert
^3 {k_1}(u)\,\mathrm{d}u
+ I(f) \int u^2 {k_1}(u)\,\mathrm{d}u \biggr\}
=0.
\end{eqnarray*}
This
ensures the validity of (\ref{Zwet2}),
and completes the proof of the contiguity in present case.

Next, consider the case where $\mathbf x_{n1},\ldots,\mathbf x_{nn}$
are nonrandom, and we observe $\mathbf
w_{ni}=\mathbf
x_{ni}+\mathbf v_{ni}, i=1,\ldots,n$.
Let ${k_2}$ denote the density of $(\mathbf v_{ni}-\bar{\mathbf
v}_n)^{\top}\mathbf b^0, i=1,\ldots,n$.
Again, by (\ref{Zwet4}),
\begin{eqnarray*}
&& \sum_{i=1}^n
H^2(P_{ni},Q_{ni})
\\
&&\quad \leq \sum_{i=1}^n\int\biggl\{\int
\bigl[f\bigl(e-n^{-1/2}u\bigr)-f(e) \bigr]{k_2}\bigl(u+(
\mathbf x_{ni}-\bar{\mathbf x}_n)^{\top}\mathbf
b^{0*}\bigr)\,\mathrm{d}u \biggr\}^2 \frac{\mathrm{d}e}{f(e)}
\nonumber
\\
&&\quad \leq \sum_{i=1}^n\int\biggl\{\int
\bigl[f\bigl(e-n^{-1/2}u\bigr) - f(e) \bigr]^2
{k_2} \bigl(u-(\mathbf x_{ni}-\bar{\mathbf
x}_n)^{\top}\mathbf b_0^* \bigr)\,\mathrm{d}u \biggr\}
\frac{\mathrm{d}e}{f(e)}
\\
&&\quad \le 2n^{1/2} \int\!\!\int_{n^{-1/2}}^{n^{-1/2}} \int
\frac
{\llvert f'(e-tu)\rrvert ^2}{f(e)} \,\mathrm{d}e \,\mathrm{d}t\, n^{-1}\sum_{i=1}^n
u^2 {k_2} \bigl(u-(\mathbf x_{ni}-\bar{\mathbf
x}_n)^{\top}\mathbf b^{0*} \bigr)\,\mathrm{d}u.
\nonumber
\end{eqnarray*}
Hence, by \textup{(F.2)} and by the change of variable formula,
%
%e3.8 #&#
%
\begin{eqnarray}
\label{Zwet5} \limsup_n \sum_{i=1}^n
H^2(P_{ni},Q_{ni}) \leq C \biggl[\int
u^2 {k_2}(u)\,\mathrm{d}u+\mathbf b^{0*\top}\mathbf
Q_n\mathbf b^{0*} \biggr]<\infty.
\end{eqnarray}
Similarly one verifies (\ref{Zwet2}) here. %\hfill$\blacksquare$
\end{pf}

Lemmas \ref{Lemma1} and \ref{Lemma2} enable us to extend the
approximation of the rank statistic $\mathbf S_n(\mathbf
b_n^*,\bolds\beta_n)$ by a sum of independent r.v.'s under
the contiguous sequence of distributions. Let
\begin{eqnarray*}
\mathbf T_n\bigl(\mathbf b_n^*,\bolds\beta_n
\bigr)&=&n^{-1/2}\sum_{i=1}^n(
\mathbf w_{ni}-\bar{\mathbf w}_n)\varphi\bigl(F
\bigl(e_{ni}-(\mathbf w_{ni}-\bar{\mathbf w}_n)^{\top
}
\mathbf b_n^*-(\mathbf v_{ni}-\bar{\mathbf
v}_n)^{\top}\bolds\beta_n\bigr) \bigr).
\end{eqnarray*}
We have the following
corollary.

\begin{corollary}\label{Cor1} Under the conditions of Theorem \ref{Theorem0},
and under $\{Q_n^{(n)}\}$,
%
%e3.9 #&#
%
\begin{eqnarray}
\label{21} \mathbf S_n\bigl(\mathbf b_n^*,\bolds
\beta_n\bigr)&=&n^{-1/2}\sum_{i=1}^n(
\mathbf w_{ni}-\bar{\mathbf w}_n)a_n \bigl(R
\bigl(e_{ni}-(\mathbf w_{ni}-\bar{\mathbf w}_n)^{\top
}
\mathbf b_n^*-(\mathbf v_{ni}-\bar{\mathbf
v}_n)^{\top}\bolds\beta_n\bigr) \bigr)
\nonumber\\[-8pt]\\[-8pt]\nonumber
&=&\mathbf T_n\bigl(\mathbf b_n^*,\bolds
\beta_n\bigr)+\mathrm{o}_p(1).
\end{eqnarray}
Hence,
\[
\mathbf S_n\bigl(\mathbf b_n^*,\bolds\beta_n
\bigr)-\mathbf S_n(\mathbf0,\mathbf0)=\mathbf T_n\bigl(
\mathbf b_n^*,\bolds\beta_n\bigr)-\mathbf T_n(
\mathbf0,\mathbf0)+\mathrm{o}_p(1).
\]
\end{corollary}

%s3.3 #&#
\subsection{Asymptotic linearity of \texorpdfstring{$\mathbf{S}_n(\mathbf{b}_n^*,\bolds{\beta}_n)$}{{S}n({b}n^*,betan)}}\label{sec3.3}

\begin{lemma}\label{Lemma3}
Under the conditions of Theorem \ref{Theorem0},
%
%e3.10 #&#
%
\begin{equation}
\label{1969} \bigl\llVert\mathbf S_n\bigl(\mathbf
b_n^*,\bolds\beta_n\bigr)-\mathbf S_n(\mathbf
0,\mathbf0) + \gamma\bigl[(\mathbf Q+\mathbf V)\mathbf b^{0*}+\mathbf V
\bolds\beta^0\bigr]\bigr\rrVert\stackrel{p} {\rightarrow} 0,
\end{equation}
where
%
%e3.11 #&#
%
\begin{equation}
\label{gamma} \gamma=\int_0^1-
\frac{f^{\prime}(F^{-1}(u))}{f(F^{-1}(u))}\varphi(u)\,\mathrm{d}u =-\int_{\mathbb
R^1}\varphi\bigl(F(z)
\bigr) \,\mathrm{d}f(z).
\end{equation}
\end{lemma}

\begin{pf}
Consider the sequence of functions $ \{\varphi
^{(k)}(\cdot) \}_{k=1}^{\infty}$
%
%e3.12 #&#
%
\begin{equation}
\label{k} \varphi^{(k)}(u)=\varphi\biggl(\frac{1}{k+1}
\biggr)\mathbb I \biggl[u<\frac{1}k \biggr]+\varphi(u)\mathbb I \biggl[
\frac{i-1}{k+1} <u\leq\frac{i}{k+1} \biggr], \qquad i=2,\ldots,k.
\end{equation}
Then, by Lemma V.1.6.a \cite{Hajek1967}, $\varphi^{(k)}$ is
nondecreasing and bounded on $(0,1)$ and
%
%e3.13 #&#
%
\begin{equation}
\label{k1} \lim_{n\rightarrow\infty}\int_0^1
\bigl[\varphi^{(k)}(u)-\varphi(u)\bigr]^2\,\mathrm{d}u=0.
\end{equation}
The function $\varphi^{(k)}$ has at most countable set $B_k$ of
discontinuity points.
Observe that assumption~\textup{(V.3)} implies that $n^{-1/2}\max_{1\le
i\le n}\{\llVert \mathbf w_{ni} -\bar{\mathbf w}_n\rrVert + \llVert
\mathbf v_{ni} -\bar
{\mathbf v}_n\rrVert \}\stackrel{p}{\rightarrow} 0$.
This fact together with the uniform continuity of $F$ implies that
\[
\sup_{e\in\R, 1\le i\le n}\bigl\llvert F\bigl(e-n^{-1/2}(\mathbf
w_{ni}-\bar{\mathbf w}_n)^{\top}\mathbf
b^{0*}-n^{-1/2}(\mathbf v_{ni}-\bar{\mathbf
v}_n)^{\top}\bolds\beta^0\bigr)- F(e)\bigr\rrvert
\stackrel{p} {\rightarrow} 0. %
\]
Hence,
\[
\varphi^{(k)} \bigl(F\bigl(e-n^{-1/2}(\mathbf w_{ni}-\bar{
\mathbf w}_n)^{\top}\mathbf b^{0*}-n^{-1/2}(
\mathbf v_{ni}-\bar{\mathbf v}_n)^{\top}\bolds
\beta^0 \bigr)\bigr)
\]
{converges to} $\varphi^{(k)} (F(e) )$, {in probability,}
uniformly in $i=1,\ldots,n$. It, {in turn}, implies that the
conditional expectation
\begin{eqnarray*}
\mathbb E \bigl[ \bigl(\varphi^{(k)} \bigl(F\bigl(e_{ni}-n^{-1/2}(
\mathbf w_{ni}-\bar{\mathbf w}_n)^{\top}\mathbf
b^{0*}-n^{-1/2}(\mathbf v_{ni}-\bar{\mathbf
v}_n)^{\top}\bolds\beta^0\bigr) \bigr)-
\varphi^{(k)} \bigl(F(e_{ni}) \bigr) \bigr)^2\mid
\mathbf v_{ni},\mathbf x_{ni} \bigr]
\end{eqnarray*}
converges to 0, {in probability,}
uniformly in $i=1,\ldots,n$ and $k$.

Let $\mathbf S_n^{(k)}(\mathbf b^*,\bolds\beta)$ and $\mathbf
T_n^{(k)}(\mathbf b^*,\bolds\beta)$ be analogous
to
$\mathbf S_n(\mathbf b^*,\bolds\beta), \mathbf T_n(\mathbf
b^*,\bolds\beta)$,
respectively, with $\varphi$ replaced with $\varphi^{(k)}$. Then\vspace*{1pt} we
can bound the norm of the covariance matrix of
$\mathbf T_n^{(k)}(\mathbf b^*,\bolds\beta)-\mathbf T_n^{(k)}(\mathbf
0,\mathbf0)$ for any fixed $k$ in the
following way:
denote
\[
\mathbf A_n^{(k)}=\mathbb E \bigl\{\bigl[\mathbf
T_n^{(k)}\bigl(\mathbf b_n^*,\bolds
\beta_n\bigr)-\mathbf T_n^{(k)}(\mathbf0,
\mathbf0)\bigr] \bigl[\mathbf T_n^{(k)}\bigl(\mathbf
b_n^*,\bolds\beta_n\bigr)-\mathbf T_n^{(k)}(
\mathbf0,\mathbf0)\bigr]^{\top} \bigr\}.
\]
Then
%
%e3.14 #&#
%
\begin{eqnarray}
\label{A} \mathbf A_n^{(k)}&=&\mathbb E \Biggl
\{n^{-1}\sum_{i=1}^n (\mathbf
w_{ni}-\bar{\mathbf w}_n) (\mathbf w_{ni}-\bar{
\mathbf w}_n)^{\top}\nonumber
\\
&&\hspace*{43pt}{}\times \bigl[\varphi^{(k)} \bigl(F\bigl(e_{ni}-(\mathbf
w_{ni}-\bar{\mathbf w}_n)^{\top}\mathbf
b_n^*-(\mathbf v_{ni}-\bar{\mathbf v}_n)^{\top}
\bolds\beta_n\bigr) \bigr)-\varphi^{(k)} \bigl(F(e_{ni})
\bigr) \bigr]^2 \Biggr\}\nonumber
\\
&=&n^{-1}\sum_{i=1}^n\mathbb E
\bigl\{(\mathbf w_{ni}-\bar{\mathbf w}_n) (\mathbf
w_{ni}-\bar{\mathbf w}_n)^{\top}
\\
&&{}\hspace*{43pt}{}\times \mathbb E \bigl[ \bigl(\varphi^{(k)} \bigl(F
\bigl(e_{ni}-(\mathbf w_{ni}-\bar{\mathbf w}_n)^{\top}
\mathbf b_n^*-(\mathbf v_{ni}-\bar{\mathbf
v}_n)^{\top}\bolds\beta_n\bigr) \bigr)\nonumber
\\
&&\hspace*{203pt}{} -\varphi
^{(k)} \bigl(F(e_{ni}) \bigr) \bigr)^2\mid\mathbf
v_{ni},\mathbf x_{ni} \bigr] \bigr\}.
\nonumber
\end{eqnarray}
Hence,
\begin{eqnarray*}
\bigl\llVert\mathbf A_n^{(k)}\bigr\rrVert&\leq&\Biggl\{
\Biggl\llVert n^{-1}\sum_{i=1}^n
(\mathbf w_{ni}-\bar{\mathbf w}_n) (\mathbf
w_{ni}-\bar{\mathbf w}_n)^{\top}-(\mathbf Q+\mathbf
V)\Biggr\rrVert+\llVert\mathbf Q+\mathbf V\rrVert\Biggr\}\cdot \mathrm{o}(1)
\\
&=& \bigl\{\llVert\mathbf Q+\mathbf V\rrVert+\mathrm{o}(1) \bigr\}\cdot \mathrm{o}(1).
\end{eqnarray*}

This{, together with the fact $\mathbb E\mathbf T_n^{(k)}(\mathbf
0,\mathbf0)=\mathbf0$,} implies
%
%e3.15 #&#
%
\begin{equation}
\label{T} \bigl\llVert\mathbf T_n^{(k)}\bigl(\mathbf
b_n^*,\bolds\beta_n\bigr)-\mathbf T_n^{(k)}(
\mathbf0,\mathbf0)-\mathbb E\mathbf T_n^{(k)}\bigl(\mathbf
b_n^*,\bolds\beta_n\bigr)\bigr\rrVert\stackrel{p} {
\rightarrow}0.
\end{equation}
Furthermore, for any fixed $k$ and for fixed $\mathbf b^{0*},\bolds
\beta^0$,
%
%e3.16 #&#
%
\begin{equation}
\label{gammak} \mathbf T_n^{(k)}\bigl(\mathbf
b_n^*,\bolds\beta_n\bigr)-\mathbf T_n^{(k)}(
\mathbf0,\mathbf0)+\gamma_k\bigl[(\mathbf Q+\mathbf
V)\mathbf
b^{0*}+\mathbf V\bolds\beta^0\bigr]\stackrel{p} {
\rightarrow} \mathbf0,
\end{equation}
where
\[
\gamma_k=-\int_{\mathbb R^1} \varphi^{(k)}
\bigl(F(e)\bigr)f^{\prime}(e)\,\mathrm{d}e=-\int_0^1
\varphi^{(k)}(u)\frac
{f^{\prime}(F^{-1}(u))}{f(F^{-1}(u))}\,\mathrm{d}u.
\]
Indeed, (we put $\bar{\mathbf x}_n=\bar{\mathbf v}_n=\mathbf0$, for
the sake of brevity)
\begin{eqnarray*}
&& n^{-1/2}\sum_{i=1}^n
\mathbb E \bigl\{\mathbf w_{ni} \bigl(\mathbb E \bigl[\varphi^{(k)}\bigl(F\bigl(e_{ni}-n^{-1/2}
\bigl(\mathbf w_{ni}^{\top}\mathbf b^{0*}-n^{-1/2}
\mathbf v_{ni}^{\top}\bolds\beta^0\bigr)-
\varphi^{k}\bigl(F(e_{ni})\bigr)
\\[-2pt]
&&\hspace*{151pt}{} -\gamma_k\bigl(n^{-1/2}\bigl[\mathbf
w_{ni}^{\top}\mathbf b^{0*}+\mathbf
v_{ni}^{\top}\bolds\beta^0\bigr]\bigl)\mid\mathbf
v_{ni},\mathbf x_{ni} \bigr] \bigr) \bigr\}
\\[-2pt]
&&\quad =n^{-1/2}\sum_{i=1}^n\mathbb E
\biggl\{\mathbf w_{ni} \biggl(\int_{\mathbb R^1}
\varphi^{(k)}\bigl(F(z)\bigr)\,\mathrm{d}\bigl[F\bigl(z+n^{-1/2}\mathbf
w_{ni}^{\top}\mathbf b^{0*}+n^{-1/2}\mathbf
v_{ni}^{\top}\bolds\beta^0\bigr)-F(z)\bigr]
\\[-2pt]
&&\hspace*{120pt}\qquad{}-n^{-1/2}\bigl[\mathbf w_{ni}^{\top}\mathbf
b^{0*}+\mathbf v_{ni}^{\top
}\bolds\beta^0
\bigr]\int_{\mathbb
R^1}\varphi^{(k)}\bigl(F(z)
\bigr)f^{\prime}(z)\,\mathrm{d}z \biggr) \biggr\}
\\[-2pt]
&&\quad =n^{-1/2}\sum_{i=1}^n\mathbb E
\biggl\{\mathbf w_{ni} \biggl(\int_{\mathbb R^1}
\varphi^{(k)}\bigl(F(z)\bigr)
\,\mathrm{d} \bigl[F\bigl(z+n^{-1/2}\mathbf w_{ni}^{\top}
\mathbf b^{0*}+n^{-1/2}\mathbf v_{ni}^{\top}
\bolds\beta^0\bigr)-F(z)
\\[-2pt]
&&\hspace*{208pt}{}-n^{-1/2}\bigl(\mathbf
w_{ni}^{\top
}\mathbf b^{0*}+\mathbf
v_{ni}^{\top}\bolds\beta^0\bigr)f(z) \bigr] \biggr)
\biggr\}
\rightarrow\mathbf0.
\end{eqnarray*}
Moreover, we have
%
%e3.17 #&#
%
\begin{eqnarray}
\label{gamma2} (\gamma_k-\gamma)^2&=& \biggl\langle(
\varphi_k-\varphi),-\frac
{f^{\prime}(F^{-1}(\cdot))}{f(F^{-1}(\cdot))} \biggr\rangle^2\nonumber
\\
&\leq&\llVert\varphi_k-\varphi\rrVert^2 \biggl\llVert
-\frac{f^{\prime
}(F^{-1}(\cdot))}{f(F^{-1}(\cdot))}\biggr\rrVert^2
\\
&=& \mathcal I(f)\llVert
\varphi_k-\varphi\rrVert^2\rightarrow0\qquad\mbox{as } k
\rightarrow\infty.
\nonumber
\end{eqnarray}
Using (\ref{gammak}), (\ref{gamma2}), Lemmas~\ref{Lemma1} and~\ref{Lemma2},
Corollary \ref{Cor1} and Lemma 3.5 in \cite{Jur1969}, we obtain that
\[
P \bigl(\bigl\llVert\mathbf S_n\bigl(\mathbf b_n^*,
\bolds\beta_n\bigr)-\mathbf S_n^{(k)}\bigl(
\mathbf b_n^*,\bolds\beta_n\bigr)\bigr\rrVert>\varepsilon
\bigr)<\varepsilon
\]
for $\forall k>k(\varepsilon)$, $\forall n>n(k)$, and finally we arrive
at (\ref{1969}). %\hfill$\blacksquare$
\end{pf}

%s3.4 #&#
\subsection{Uniform asymptotic quadraticity of the Jaeckel dispersion}\label{sec3.4}
Recall that $\bar a_n=0$ under \textup{(A.1)}. Rewrite the Jaeckel
dispersion in the presence of measurement errors in the form
\begin{eqnarray*}
\label{S6a} \mathcal D_n(\mathbf{b})&=& \sum
_{i=1}^n \bigl[Y_{ni}-
\mathbf{w}_{ni}^{\top}\mathbf{b} \bigr]a_n
\bigl({R}_{ni}(\mathbf{b})\bigr)
\end{eqnarray*}
or eventually in the form
\begin{eqnarray*}
&&\mathcal D_n\bigl(\mathbf{b^*,\bolds\beta}\bigr)=\sum
_{i=1}^n \bigl[e_{ni}-(
\mathbf{w}_{ni}-\bar{\mathbf w}_n)^{\top}
\mathbf{b}^*-(\mathbf{v}_{ni}-\bar{\mathbf v}_n)^{\top
}
\bolds\beta\bigr]a_n \bigl(R\bigl(e_{ni}-\mathbf
w_{ni}^{\top}\mathbf{b^*}-\mathbf v_{ni}^{\top}
\bolds\beta\bigr)\bigr),
\nonumber
\end{eqnarray*}
where $\mathbf b^*=\mathbf b-\bolds\beta$.
It is a piecewise linear, convex function of $\mathbf b$ and $\mathbf
b^*$, respectively. Hence, the minimum
$\widehat{\bolds\beta}_n=\arg\min_{\mathbf b\in\mathbb
R^p}\mathcal D_n(\mathbf{b})$ exists, and is considered as an estimate of
$\bolds\beta$ in model (\ref{2}). By \cite{Jaeckel}, the partial
derivatives of $\mathcal D_n(\mathbf b)$ exist for almost all $\mathbf
b$, and where they exist, are equal to
\[
\frac{\partial}{\partial b_j}\mathcal D_n(\mathbf b)=-n^{1/2}S_{nj}(
\mathbf b)=-\sum_{i=1}^n(w_{nij}-
\bar{w}_j)a_n(Y_i-\mathbf w_{ni}
\mathbf b),\qquad j=1,\ldots,p.
\]
Otherwise speaking,
\begin{eqnarray*}
&&\nabla\mathcal D_n\bigl(\mathbf b^*,\bolds\beta
\bigr)=-n^{1/2}\mathbf S_{n}\bigl(\mathbf b^*,\bolds\beta
\bigr)=-\sum_{i=1}^n(\mathbf
w_{ni}-\bar{\mathbf w}_n)a_n\bigl(R
\bigl(e_{ni}-\mathbf w_{ni}^{\top}\mathbf b^*-\mathbf
v_{ni}^{\top}\bolds\beta\bigr)\bigr),
\end{eqnarray*}
where $\nabla$ denotes the subgradient.

Consider the quadratic function
\begin{eqnarray*}
\label{quadratic} \mathcal C_n\bigl(\mathbf b^*,\bolds\beta\bigr)=
\tfrac{1}{2}\gamma\mathbf b^{*\top}(\mathbf Q+\mathbf V)\mathbf b^*-
\mathbf b^{*\top}\mathbf S_{n}(\mathbf0)+\gamma\mathbf b^*
\mathbf V\bolds\beta+\mathcal D_n(\mathbf0).
\end{eqnarray*}
Then $\mathcal D_n(\mathbf b)$ and $\mathcal C_n(\mathbf b)$ are both
convex functions and $\mathcal D_n(\mathbf0)=\mathcal C_n(\mathbf0)$.
Moreover,
\begin{eqnarray*}
\label{41} \nabla\bigl[\mathbf D_n\bigl(\mathbf b^*,\bolds\beta
\bigr)-\mathcal C_n\bigl(\mathbf b^*,\bolds\beta\bigr)\bigr]=-
\bigl[n^{1/2}\bigl(\mathbf S_n\bigl(\bigl(\mathbf
b^*,\bolds\beta
\bigr)\bigr)-\mathbf S_n(\mathbf0,\mathbf0)+\gamma(\mathbf Q+\mathbf
V)\mathbf b^*+\gamma\mathbf V\bolds\beta\bigr)\bigr].
\end{eqnarray*}
Hence, it follows from (\ref{1969}) that for $\mathbf
b_n^*=n^{-1/2}\mathbf b^{0*}, \bolds\beta_n=n^{-1/2}\bolds\beta
^0$ with $\mathbf b^{0*}, \bolds\beta^0\in\mathbb R^p$ fixed that
\[
\bigl\llVert\nabla\bigl[\mathbf D_n\bigl(n^{-1/2}\mathbf
b^{0*},n^{-1/2}\bolds\beta^0\bigr)-\mathcal
C_n\bigl(n^{-1/2}\mathbf b^{0*},n^{-1/2}
\bolds\beta^0\bigr)\bigr]\bigr\rrVert\stackrel{p} {\rightarrow}0.
\]
Using the convexity arguments in \cite{Heiler} (Appendix) and \cite{Pollard91} (Convexity lemma), we conclude that
\[
\sup\bigl\llvert\mathcal D_n\bigl(n^{-1/2}\mathbf
b^{0*},n^{-1/2}\bolds\beta^0\bigr)-
\tfrac{1}{2}\gamma\mathbf b^{0*\top}(\mathbf Q+\mathbf V)\mathbf
b^{0*} +\mathbf b^{0*\top}\mathbf S_{n}(\mathbf0)-
\gamma\mathbf b^{0*}\mathbf V\bolds\beta^0+\mathcal
D_n(\mathbf0)\bigr\rrvert=\mathrm{o}_p(1), %
\]
{where the supremum is taken over the set $\{{\llVert \mathbf
b^{0*}\rrVert \leq C,
\llVert \bolds\beta^0\rrVert \leq C}\}$.}
Hence, following the arguments in the proof of Theorem 1 in \cite
{Pollard91}, we conclude that, under the local alternative $\bolds
\beta_n=n^{-1/2}\bolds\beta^0$,
\begin{eqnarray*}
\label{43} \mathop{\arg\min}_{\mathbf b^{0*}}\mathcal D_n
\bigl(n^{-1/2}\mathbf b^{0*},n^{-1/2}\bolds
\beta^0\bigr)
\end{eqnarray*}
is asymptotically equivalent to
%
%e3.18 #&#
%
\begin{equation}
\label{44} \mathop{\arg\min}_{\mathbf b^{0*}} \biggl[\frac{1}{2}\gamma
\mathbf b^{0*\top
}(\mathbf Q+\mathbf V)\mathbf b^{0*}-\mathbf
b^{0*\top}\mathbf S_{n}(\mathbf0) +\gamma\mathbf
b^{0*}\mathbf V\bolds\beta^0 \biggr].
\end{equation}
The solution of (\ref{44}) equals to
\[
\mathbf b^{0*}=\mathbf b^0-\bolds\beta^0=n^{1/2}(
\widehat{\bolds\beta}_n-\bolds\beta_n)=
\gamma^{-1}(\mathbf Q+\mathbf V)^{{-1}}\mathbf S_n(
\mathbf0,\mathbf0)-(\mathbf Q+\mathbf V)^{{-1}}\mathbf V\bolds
\beta^0.
\]
Hence, in the linear model with local value of regression parameter
$\beta$, when
\[
Y_{ni}=\mathbf x_{ni}^{\top}\bolds
\beta_n+e_{ni},\qquad \bolds\beta_n=n^{-1/2}
\bolds\beta^0,
\]
when we observe only $\mathbf w_{ni}=\mathbf x_{ni}+\mathbf v_{ni}$
instead of $\mathbf x_{ni}, i=1,\ldots,n$,
the R-estimator is asymptotically normally distributed with a bias
$\mathbf B=-(\mathbf Q+\mathbf V)^{-1}\mathbf V\bolds\beta^0$, that is,
%
%e3.19 #&#
%
\begin{equation}
\label{normality1} n^{1/2}\bigl(\widehat{\bolds\beta}_n-n^{-1/2}
\bolds\beta^0\bigr)\stackrel{\mathcal D} {\rightarrow}\mathcal
N_p \bigl(B,(\mathbf Q+\mathbf V)^{-1} A^2(
\varphi) \bigr),\qquad \mathbf B= - (\mathbf Q+\mathbf V)^{-1}\mathbf V\bolds
\beta^0.
\end{equation}

Finally, as we have already mentioned,
all of the above arguments and motivations are valid when we replace
$e_{ni}$ with $e_{ni}+u_{ni}, i=1,\ldots,n$. This completes the
proof of Theorem \ref{Theorem0}. %\hfill$\blacksquare$
%\end{pf}

%s4 #&#
\section{Numerical illustration}\label{sec4}
The following simulation study illustrates the effect of measurement
errors in regressors on the finite-sample performance of R-estimates.
Empirical bias (and variance) of R-estimates are computed and compared
for various measurement error models. For the sake of comparison, the
bias{es} and variance{s} are {also} computed for the least squares
estimate (LSE) and {the least absolute deviation
($L_1$) estimate,} under the same setup. {Moreover, we compare the
deterministic and random regressors.}

All the simulations were performed in the statistical software \verb R  using
standard tools and libraries. For minimization of (\ref{S6})
functions \verb optimize ~and \verb optim ~with initial estimate 0.5
-- regression quantile were used.
The random numbers generator was setup with the initial value \verb
set.seed(15).

The results illustrate that the bias of R-estimate is surprisingly
stable with respect to the sample size; the bias corresponding to small
$n$ is comparable to the asymptotic one derived in Theorem~\ref{Theorem0}.

Notice that the bias of R-estimator
only slightly differs from the biases of LSE and \mbox{$L_1$-}estimators.

%s4.1 #&#
\subsection{Regression line}\label{sec4.1}

Consider first the model of regression line
\[
Y_{i}=\beta_0+x_{i}\beta_{1}+e_{i},
\qquad i=1,\ldots,n,
\]
where {the $Y_i$ are measured accurately}, while instead of $x_i$ we
observe only $w_i=x_i+v_i, i=1,\ldots,n$.
The R-estimator of parameter $\beta_{1}$ is based on
Wilcoxon scores generated by score function $\varphi(u)=u-1/2$.

%t1 #&#
%
\begin{table}
\caption{Empirical bias of R-estimator for various $n$ and
measurement errors $v_{i}$; nonrandom regressors $x_{i}$}\label{T1}
\begin{tabular*}{\tablewidth}{@{\extracolsep{\fill}}@{}lllllllll@{}}
\hline
& \multicolumn{8}{l@{}}{$n$}\\[-8pt]
& \multicolumn{8}{l@{}}{\hrulefill}\\
$v_{i}$ & 10 & 20 & 50 & 100 & 200 & 500 & 1000 & $\infty$\\
\hline
0& \phantom{$-$}0.002& \phantom{$-$}0.001&\phantom{$-$}0.000&\phantom{$-$}0.000&\phantom{$-$}0.000&\phantom{$-$}0.000&\phantom{$-$}0.000&\phantom{$-$}0.000\\
$\mathcal{U}(-5,0)$&$-$0.264&$-$0.295&$-$0.305&$-$0.297&$-$0.302&$-$0.306&$-$0.307&$-$0.296\\
$\mathcal{U}(0,9)$&$-$0.684&$-$0.727&$-$0.732&$-$0.714&$-$0.719&$-$0.727&$-$0.728&$-$0.720\\
$\mathcal{U}(-3,9)$&$-$0.982&$-$1.013&$-$1.006&$-$0.983&$-$0.986&$-$0.995&$-$0.995&$-$1.000\\
$\mathcal{N}(0,1)$&$-$0.128&$-$0.148&$-$0.150&$-$0.146&$-$0.148&$-$0.151&$-$0.152&$-$0.154\\
$\mathcal{N}(0,2)$&$-$0.440&$-$0.483&$-$0.488&$-$0.476&$-$0.480&$-$0.487&$-$0.488&$-$0.500\\
$\mathcal{N}(0,3)$&$-$0.790&$-$0.836&$-$0.837&$-$0.819&$-$0.822&$-$0.832&$-$0.833&$-$0.857\\
\hline
\end{tabular*}
\end{table}

%t2 #&#
%
\begin{table}[b]
\tabcolsep=0pt
\caption{Empirical bias of R-estimator for various $n$ and
measurement errors $v_{i}$; random regressors $x_{i}$}\label{T2}
\begin{tabular*}{\tablewidth}{@{\extracolsep{\fill}}@{}lllllllll@{}}
\hline
& \multicolumn{8}{l@{}}{$n$}\\[-8pt]
& \multicolumn{8}{l@{}}{\hrulefill}\\
$v_{i}$ & 10 & 20 & 50 & 100 & 200 & 500 & 1000 & $\infty$\\
\hline
0&\phantom{$-$}0.004&$-$0.001&\phantom{$-$}0.000&\phantom{$-$}0.000&\phantom{$-$}0.000&\phantom{$-$}0.000&\phantom{$-$}0.000&\phantom{$-$}0.000\\
$\mathcal{U}(-5,0)$&$-$0.283&$-$0.297&$-$0.305&$-$0.306&$-$0.307&$-$0.309&$-$0.309&$-$0.296\\
$\mathcal{U}(0,9)$&$-$0.711&$-$0.722&$-$0.728&$-$0.730&$-$0.730&$-$0.732&$-$0.732&$-$0.720\\
$\mathcal{U}(-3,9)$&$-$0.998&$-$1.000&$-$0.999&$-$1.000&$-$0.999&$-$1.000&$-$1.000&$-$1.000\\
$\mathcal{N}(0,1)$&$-$0.138&$-$0.149&$-$0.150&$-$0.153&$-$0.153&$-$0.153&$-$0.153&$-$0.154\\
$\mathcal{N}(0,2)$&$-$0.462&$-$0.481&$-$0.487&$-$0.489&$-$0.491&$-$0.492&$-$0.492&$-$0.500\\
$\mathcal{N}(0,3)$&$-$0.813&$-$0.830&$-$0.833&$-$0.835&$-$0.837&$-$0.837&$-$0.838&$-$0.857\\
\hline
\end{tabular*}
\end{table}

All the simulation results are based on 10\,000 replications,
parameters were chosen as
$\beta_0=1, \beta_{1}=2$, and model errors $e_{i}$ follow the
logistic distribution.
{In Tables~\ref{T1} and \ref{T2}, the empirical bias of R-estimator
based on Wilcoxon scores is compared for various sample sizes ($n=10,
\ldots, 1000$) and with the
theoretical asymptotic result ($n=\infty$). The regressors $x_i$ are
deterministic in Table~\ref{T1}; they were generated from uniform
$\mathcal{U}(-3,9)$ distribution once for all experiment and then
considered as fixed. The regressors in Table~\ref{T2} are random; each
time they were generated also from uniform distribution $\mathcal
{U}(-3,9)$. This enables to see the difference between deterministic
and random regressors: the bias differs more from its asymptotic value
in case of deterministic regressors than in case of random regressors;
it can be caused by the slower rate of convergence.}
The measurement errors $v_i$ are either uniformly or normally
distributed ($i=1,\ldots,n$).

%t3 #&#
%
\begin{table}
\tabcolsep=0pt
\caption{Empirical bias (variance) of R-estimator, LSE and
$L_1$-estimator for various measurement errors $v_{i}$ and model errors
$e_{i}$; $n=50$}\label{T3}
\begin{tabular*}{\tablewidth}{@{\extracolsep{\fill}}@{}llllll@{}}
\hline
& \multicolumn{5}{l@{}}{$e_{i}$}\\[-8pt]
& \multicolumn{5}{l@{}}{\hrulefill}\\
$v_{i}$ &  Normal & Logistic & Laplace & Pareto & Cauchy\\
\hline
&\phantom{$-$}0.002 {(0.182)} &\phantom{$-$}0.008 {(0.527)}&$-$0.002 {(0.254)}&\phantom{$-$}0.000 {(0.416)}& \phantom{$-$}0.018 {(0.672)}\\
0&\phantom{$-$}0.004 {(0.172)}&\phantom{$-$}0.009 {(0.567)}&$-$0.003 {(0.355)}&\phantom{$-$}5.047 {(89\,200)}& $-$3.289 {(85\,700)}\\
&\phantom{$-$}0.002 {(0.274)} &\phantom{$-$}0.010 {(0.708)}&$-$0.002 {(0.249)}&\phantom{$-$}0.000 {(1.191)}& \phantom{$-$}0.018 {(0.526)}
\\[3pt]
&$-$0.399 {(0.155)}&$-$0.398 {(0.438)}&$-$0.396 {(0.214)}&$-$0.401 {(0.422)}&$-$0.404 {(0.568)}\\
$\mathcal{U}(-3,3)$&$-$0.395 {(0.147)}&$-$0.396 {(0.466)}&$-$0.394 {(0.283)}&$-$7.137 {(604\,000)}&\phantom{.}22.62 {(4\,000\,000)}\\
&$-$0.401 {(0.232)}&$-$0.400 {(0.591)}&$-$0.400 {(0.235)}&$-$0.422 {(0.932)}&$-$0.405 {(0.456)}
\\[3pt]
&$-$0.995 {(0.101)}&$-$1.006 {(0.278)}&$-$0.997 {(0.142)}&$-$1.001 {(0.309)}&$-$1.010 {(0.397)}\\
$\mathcal{U}(-6,6)$&$-$0.995 {(0.096)}&$-$1.009 {(0.294)}&$-$0.998 {(0.182)}&$-$7.259 {(401\,000)}& \phantom{$-$}0.933 {(36\,400)}\\
&$-$0.995 {(0.151)}&$-$1.006 {(0.376)}&$-$0.995 {(0.157)}&$-$1.001 {(0.587)}&$-$1.014 {(0.320)}
\\[3pt]
&$-$0.153 {(0.174)}&$-$0.161 {(0.493)}&$-$0.145 {(0.243)}&$-$0.149 {(0.439)}&$-$0.147 {(0.638)}\\
$\mathcal{N}(0,1)$ &$-$0.152 {(0.163)}&$-$0.158 {(0.523)}&$-$0.145 {(0.328)}&$-$2.136 {(281\,000)}&$-$7.380 {(739\,000)}\\
&$-$0.153 {(0.261)}&$-$0.159 {(0.675)}&$-$0.147 {(0.259)}&$-$0.165 {(1.092)}&$-$0.138 {(0.510)}\\
\hline
\end{tabular*}
\end{table}

{Table~\ref{T3} compares empirical bias and variance (in parenthesis)
of R-estimator based on Wilcoxon scores, of LSE and $L_1$-estimate
{for the} sample size $n=50$ and
{when}
regressors $x_{i}$ {are random}, generated from uniform $\mathcal
{U}(-3,9)$ distribution; model errors $e_{i}$ are generated from
normal, logistic, Laplace, Pareto with parameter $\alpha=0.9$ and
Cauchy distributions. The measurement errors $v_i$ follow various
distributions, similarly as in Tables~\ref{T1} and \ref{T2}.}

%s4.2 #&#
\subsection{Model of two regressors}\label{sec4.2}
Consider the model
\[
Y_{i}=\beta_0+x_{i,1}\beta_{1}+x_{i,2}
\beta_{2}+e_{i},\qquad i=1,\ldots,n,
\]
where again the $Y_i$ are measured accurately, but instead of $\mathbf
x_i$ we observe only $\mathbf w_i=\mathbf x_i+\mathbf v_i, i=1,\ldots,n$.
The R-estimator of parameter $\bolds\beta= (\beta_{1},\beta
_{2})^{\top}$ is based on
Wilcoxon scores generated by score function $\varphi(u)=u-1/2$.

%
%t4 #&#
%
\begin{table}
\tabcolsep=0pt
\caption{Empirical bias (variance) of R-estimator, LSE and $L_1$-estimator
for various measurement errors $\mathbf v_{i}$ and model errors
$e_{i}$; $n=50$}\label{T4}
\begin{tabular*}{\tablewidth}{@{\extracolsep{\fill}}@{}lllllll@{}}
\hline
& &\multicolumn{5}{l@{}}{$e_{i}$}\\[-8pt]
& &\multicolumn{5}{l@{}}{\hrulefill}\\
$v_{i}$ &  & Normal & Logistic & Laplace & Pareto & Cauchy\\
\hline
$\mathbf0$ & $\widehat{\beta}_{1}$& \phantom{$-$}0.000 {(0.600)}&$-$0.015 {(1.688)}&\phantom{$-$}0.002 {(0.821)}&\phantom{0}$-$0.002 {(0.017)}& \phantom{$-$}0.026 {(2.326)}\\
&&$-$0.001 {(0.569)}&$-$0.019 {(1.789)}&\phantom{$-$}0.006 {(1.120)}& \phantom{$-$}21.73 {(5\,330\,000)}&$-$3.698 {(48\,700)}\\
&&$-$0.007 {(0.864)}&$-$0.008 {(2.295)}&\phantom{$-$}0.002 {(0.857)}&\phantom{2}$-$0.003 {(0.034)}& \phantom{$-$}0.035 {(1.843)}
\\[3pt]
&$\widehat{\beta}_{2}$& \phantom{$-$}0.020 {(1.176)}&\phantom{$-$}0.026 {(3.497)}& \phantom{$-$}0.001 {(1.678)}&\phantom{0}$-$0.003 {(0.033)}&\phantom{$-$}0.037 {(4.634)}\\
&&\phantom{$-$}0.014 {(1.117)}&\phantom{$-$}0.027 {(3.725)}&$-$0.001 {(2.250)}&$-$29.28 {(9\,770\,000)}&\phantom{$-$}0.695 {(66\,200)}\\
&&\phantom{$-$}0.031 {(1.744)}&\phantom{$-$}0.017 {(4.588)}&$-$0.001 {(1.758)}&\phantom{0}$-$0.005 {(0.067)}&\phantom{$-$}0.007 {(3.618)}
\\[6pt]
$\mathcal{N}_{2}(\bolds\mu, \mathbf S_{3})$&
 $\widehat{\beta}_{1}$& $-$0.362 {(0.554)}&$-$0.379 {(1.603)}&$-$0.379 {(0.798)}&\phantom{0}$-$0.372 {(0.087)}&$-$0.347 {(2.161)}\\
&&$-$0.359 {(0.528)}&$-$0.385 {(1.693)}&$-$0.376 {(1.030)}&\phantom{$-$}22.81 {(6\,270\,000)}&$-$4.118 {(48\,500)}\\
&&$-$0.365 {(0.823)}&$-$0.369 {(2.134)}&$-$0.375 {(0.881)}&\phantom{0}$-$0.370 {(0.102)}&$-$0.338 {(1.758)}
\\[3pt]
& $\widehat{\beta}_{2}$& $-$0.774 {(0.936)}&$-$0.738 {(2.662)}&$-$0.754 {(1.295)}&\phantom{0}$-$0.770 {(0.136)}&$-$0.738 {(3.618)}\\
&&$-$0.776 {(0.881)}&$-$0.734 {(2.832)}&$-$0.757 {(1.696)}&$-$26.76 {(7\,920\,000)}&$-$0.268 {(78\,500)}\\
&&$-$0.769 {(1.402)}&$-$0.754 {(3.366)}&$-$0.754 {(1.409)}&\phantom{0}$-$0.769 {(0.164)}&$-$0.752 {(2.935)}
\\[6pt]
$\mathcal{N}_{2}(\bolds\mu, \mathbf S_{2})$&
 $\widehat{\beta}_{1}$ & $-$0.643 {(0.419)}&$-$0.652 {(1.155)}&$-$0.648 {(0.579)}&\phantom{0}$-$0.649 {(0.076)}&$-$0.626 {(1.622)}\\
&&$-$0.640 {(0.399)}&$-$0.655 {(1.216)}&$-$0.653 {(0.750)}&\phantom{$-$}14.66 {(3\,070\,000)}&$-$3.987 {(47\,200)}\\
&&$-$0.647 {(0.615)}&$-$0.645 {(1.573)}&$-$0.642 {(0.657)}&\phantom{0}$-$0.650 {(0.089)}&$-$0.623 {(1.317)}
\\[3pt]
&$\widehat{\beta}_{2}$& $-$0.495 {(0.643)}&$-$0.474 {(1.730)}&$-$0.477 {(0.866)}&\phantom{0}$-$0.494 {(0.107)}&$-$0.466 {(2.426)}\\
&&$-$0.497 {(0.605)}&$-$0.469 {(1.843)}&$-$0.471 {(1.139)}&$-$17.59 {(3\,660\,000)}&$-$0.646 {(55\,900)}\\
&&$-$0.489 {(0.948)}&$-$0.477 {(2.282)}&$-$0.486 {(0.944)}&\phantom{0}$-$0.492 {(0.129)}&$-$0.478 {(1.992)}
\\[6pt]
$\mathcal{N}_{2}(\bolds\mu, \mathbf S_{1})$
&$\widehat{\beta}_{1}$ & $-$0.997 {(0.329)}&$-$1.013 {(0.879)}&$-$1.010 {(0.448)}&\phantom{0}$-$1.005 {(0.071)}&$-$0.987 {(1.264)}\\
&&$-$0.999 {(0.311)}&$-$1.015 {(0.931)}&$-$1.011 {(0.577)}&\phantom{$-$1}9.435 {(1\,260\,000)}&$-$2.474 {(40\,500)}\\
&&$-$0.994 {(0.484)}&$-$1.009 {(1.192)}&$-$1.011 {(0.509)}&\phantom{0}$-$1.005 {(0.086)}&$-$0.998 {(1.037)}
\\[3pt]
&$\widehat{\beta}_{2}$& $-$0.505 {(0.670)}&$-$0.493 {(1.797)}&$-$0.496 {(0.917)}&\phantom{0}$-$0.499 {(0.141)}&$-$0.482 {(2.489)}\\
&&$-$0.505 {(0.630)}&$-$0.486 {(1.883)}&$-$0.487 {(1.168)}&$-$21.75 {(5\,680\,000)}& \phantom{$-$}0.106 {(49\,600)}\\
&&$-$0.501 {(0.980)}&$-$0.520 {(2.374)}&$-$0.504 {(1.024)}&\phantom{0}$-$0.500 {(0.170)}&$-$0.501 {(2.064)}\\
\hline
\end{tabular*}
\end{table}

{Here we chose} $n=50$, parameters
$\beta_0=1, \beta_{1}=2, \beta_{2}=1$, random regressors $\mathbf
x_{i}=(x_{i,1},x_{i,2})^{\top}$ are generated from 2-dimensional
normal distributions $\mathcal{N}_{2}(\bolds\mu,\mathbf S_{\nu}),
\nu=1,2,3$, where $\bolds\mu=(0,1)^{\top}$ and
\[
\mathbf S_{1} = \pmatrix{ 4& 0.5
\cr
0.5& 2}, \qquad\mathbf
S_{2} = \pmatrix{2& 0.2
\cr
0.2& 2}, \qquad\mathbf S_{3} =
\pmatrix{ 1& 0.9
\cr
0.9& 1}. %
\]
Table~\ref{T4} compares empirical bias and variance (in parentheses)
of R-estimator based on Wilcoxon scores, {with those of } the LSE
and $L_1$-estimator for various {distributions of the} measurement
errors $\mathbf v_{i}$ and model errors $e_{i}$.

We have also computed R-estimates generated by other score functions,
for example, van der Waerden, median; also
another simulation design was considered -- various sample sizes $n$,
values of the parameters, distributions of regressors, measurement
errors {$v_{i}$ and $u_{i}$} and model errors.
It is of interest that the results for corresponding R-estimates are
quite similar to those presented in the previous tables.

The simulation study confirms that R-estimates in measurement error
models are biased, as well as other usual estimates. The bias is
relatively stable with respect to the sample size and to distribution
of model errors.
The R-estimates provide meaningful results as long as the $e_{i}$ have
a finite Fisher information; even under the normal errors are their
empirical variances only slightly greater than that of LSE.
The bias and other properties of R-estimates are comparable with those
of the least squares and of $L_1$ estimates unless %The problem may
%arise
the distribution of model errors $e_{i}$ is heavy, where the LSE fails.
Generally, the reduction of the bias is rather a matter of measurement
precision, of calibration and repeated measurements.

\section*{Acknowledgements}

The authors thank the Referee for his/her comments, which helped to
better understanding the text.

The work of R. Navr\'atil and J. Jure\v{c}kov\'a was partially done
during their visits to Michigan State University, hosted by the
Department of Statistics and Probability.
Research of J. Jure\v{c}kov\'a was partially supported by the Grant
GA\v{C}R 201/12/0083. Research of H.L. Koul was supported in part by
the NSF DMS Grant 1250271. Research of R. Navr\'atil was supported by
the Grant SVV-2013-267 315, and J. Picek and R. Navr\'atil were
supported by the Project Klimatext CZ.1.07/2.3.00/20.0086.

%\begin{appendix}
%\section{}
%\end{appendix}

%\section*{Acknowledgements}

%\begin{supplement}%[id=suppA]
%\slink[doi]{10.3150/00-BEJXXXXSUPP} %[doi,text={...}] - jei reikia
%suskaldyti doi
%\sfilename{BEJ000\_supp.pdf}
%\end{supplement}

% imsref loaded by linak, 2014-12-10 11:03:00
%
% imsref loaded by linak, 2014-12-15 15:50:49

\printhistory

\begin{thebibliography}{32}

%b1 ###
\bibitem{Adcock}
%
\begin{barticle}[auto:parserefs-M02]
\bauthor{\bsnm{Adcock},~\bfnm{R.~J.}\binits{R.J.}}
(\byear{1877}).
\btitle{Note on the method of least squares}.
\bjournal{The Analyst}
\bvolume{4}
\bpages{183--184}.
\end{barticle}
%

\bptok{imsref}%
\endbibitem

%b2 ###
\bibitem{Akritas}
%
\begin{barticle}[auto:parserefs-M02]
\bauthor{\bsnm{Akritas},~\bfnm{M.~G.}\binits{M.G.}} \AND
\bauthor{\bsnm{Bershady},~\bfnm{M.~A.}\binits{M.A.}}
(\byear{1996}).
\btitle{Linear regression for astronomical data with measurement errors
and intrinsic scatter}.
\bjournal{Astrophysical Journal}
\bvolume{470}
\bpages{706--728}.
\end{barticle}
%

\bptok{imsref}%
\endbibitem

%b3 ###
\bibitem{Arias2001}
%
\begin{barticle}[auto:parserefs-M02]
\bauthor{\bsnm{Arias},~\bfnm{O.}\binits{O.}},
\bauthor{\bsnm{Hallock},~\bfnm{K.~F.}\binits{K.F.}} \AND
\bauthor{\bsnm{Sosa-Escudero},~\bfnm{W.}\binits{W.}}
(\byear{2001}).
\btitle{Individual heterogeneity in the returns to schooling:
Instrumental variables quantile regression using twins data}.
\bjournal{Empirical Economics}
\bvolume{26}
\bpages{7--40}.
\end{barticle}
%

\bptok{imsref}%
\endbibitem

%b4 ###
\bibitem{Carroll2007}
%
\begin{barticle}[mr]
\bauthor{\bsnm{Carroll},~\bfnm{Raymond~J.}\binits{R.J.}},
\bauthor{\bsnm{Delaigle},~\bfnm{Aurore}\binits{A.}} \AND
\bauthor{\bsnm{Hall},~\bfnm{Peter}\binits{P.}}
(\byear{2007}).
\btitle{Non-parametric regression estimation from data contaminated by
a mixture of {B}erkson and classical errors}.
\bjournal{J. R. Stat. Soc. Ser. B Stat. Methodol.}
\bvolume{69}
\bpages{859--878}.
\bid{doi={10.1111/j.1467-9868.2007.00614.x}, issn={1369-7412}, mr={2368574}}
\end{barticle}
%

\bptok{imsref}%
% NOT OUTPUTTED:
% number = 5
% doi = http://dx.doi.org/10.1111/j.1467-9868.2007.00614.x
% fjournal = Journal of the Royal Statistical Society. Series B.
%Statistical Methodology
\endbibitem

%b5 ###
\bibitem{Carroll1999}
%
\begin{barticle}[mr]
\bauthor{\bsnm{Carroll},~\bfnm{Raymond~J.}\binits{R.J.}},
\bauthor{\bsnm{Maca},~\bfnm{Jeffrey~D.}\binits{J.D.}} \AND
\bauthor{\bsnm{Ruppert},~\bfnm{David}\binits{D.}}
(\byear{1999}).
\btitle{Nonparametric regression in the presence of measurement error}.
\bjournal{Biometrika}
\bvolume{86}
\bpages{541--554}.
\bid{doi={10.1093/biomet/86.3.541}, issn={0006-3444}, mr={1723777}}
\end{barticle}
%

\bptok{imsref}%
% NOT OUTPUTTED:
% number = 3
% doi = http://dx.doi.org/10.1093/biomet/86.3.541
% coden = BIOKAX
% fjournal = Biometrika
\endbibitem

%b6 ###
\bibitem{Carroll2006}
%
\begin{bbook}[mr]
\bauthor{\bsnm{Carroll},~\bfnm{Raymond~J.}\binits{R.J.}},
\bauthor{\bsnm{Ruppert},~\bfnm{David}\binits{D.}},
\bauthor{\bsnm{Stefanski},~\bfnm{Leonard~A.}\binits{L.A.}} \AND
\bauthor{\bsnm{Crainiceanu},~\bfnm{Ciprian~M.}\binits{C.M.}}
(\byear{2006}).
\btitle{Measurement Error in Nonlinear Models. A Modern Perspective},
\bedition{2nd} ed.
\bseries{Monographs on Statistics and Applied Probability}
\bvolume{105}.
\blocation{Boca Raton, FL}:
\bpublisher{Chapman \& Hall/CRC}.
\bid{doi={10.1201/9781420010138}, mr={2243417}}
\end{bbook}
%

\bptok{imsref}%
% NOT OUTPUTTED:
% doi = http://dx.doi.org/10.1201/9781420010138
% isbn = 978-1-58488-633-4; 1-58488-633-1
% fpage = xxviii+455
\endbibitem

%b7 ###
\bibitem{vanNess1999}
%
\begin{bbook}[mr]
\bauthor{\bsnm{Cheng},~\bfnm{Chi-Lun}\binits{C.-L.}} \AND
\bauthor{\bsnm{Van Ness},~\bfnm{John~W.}\binits{J.W.}}
(\byear{1999}).
\btitle{Statistical Regression with Measurement Error}.
\bseries{Kendall's Library of Statistics}
\bvolume{6}.
\blocation{London}:
\bpublisher{Arnold}.
\bid{mr={1719513}}
\end{bbook}
%

\bptok{imsref}%
% NOT OUTPUTTED:
% isbn = 0-340-61461-7
% fpage = xiv+262
\endbibitem

%b8 ###
\bibitem{FanTruong1993}
%
\begin{barticle}[auto]
\bauthor{\bsnm{Fan},~\bfnm{Jianqing}\binits{J.}} \AND
\bauthor{\bsnm{Truong},~\bfnm{Young~K.}\binits{Y.K.}}
(\byear{1993}).
\btitle{Nonparametric regression estimation involving errors-in-variables}.
\bjournal{Ann. Statist.}
\bvolume{21}
\bpages{23--37}.
\end{barticle}
%

\bptok{imsref}%
% NOT OUTPUTTED:
% number = 4
% doi = http://dx.doi.org/10.1214/aos/1176349402
% coden = ASTSC7
% fjournal = The Annals of Statistics
\endbibitem

%b9 ###
\bibitem{Fuller1987}
%
\begin{bbook}[mr]
\bauthor{\bsnm{Fuller},~\bfnm{Wayne~A.}\binits{W.A.}}
(\byear{1987}).
\btitle{Measurement Error Models}.
\bseries{Wiley Series in Probability and Mathematical Statistics:
Probability and Mathematical Statistics}.
\blocation{New York}:
\bpublisher{Wiley}.
\bid{doi={10.1002/9780470316665}, mr={0898653}}
\end{bbook}
%

\bptok{imsref}%
% NOT OUTPUTTED:
% doi = http://dx.doi.org/10.1002/9780470316665
% isbn = 0-471-86187-1
% fpage = xxiv+440
\endbibitem

%b10 ###
\bibitem{Hajek1967}
%
\begin{bbook}[mr]
\bauthor{\bsnm{H{\'a}jek},~\bfnm{Jaroslav}\binits{J.}} \AND
\bauthor{\bsnm{{\v{S}}id{\'a}k},~\bfnm{Zbyn{\v{e}}k}\binits{Z.}}
(\byear{1967}).
\btitle{Theory of Rank Tests}.
\blocation{New York}:
\bpublisher{Academic Press}.
\bid{mr={0229351}}
\end{bbook}
%

\bptok{imsref}%
% NOT OUTPUTTED:
% fpage = 297
\endbibitem

%b11 ###
\bibitem{Hausman2001}
%
\begin{barticle}[auto:parserefs-M02]
\bauthor{\bsnm{Hausman},~\bfnm{J.}\binits{J.}}
(\byear{2001}).
\btitle{Mismeasured variables in econometric analysis: Problems from
the right and problems from the left}.
\bjournal{J. Econ. Perspect.}
\bvolume{15}
\bpages{57--67}.
\end{barticle}
%

\bptok{imsref}%
\endbibitem

%b12 ###
\bibitem{HeLiang}
%
\begin{barticle}[mr]
\bauthor{\bsnm{He},~\bfnm{Xuming}\binits{X.}} \AND
\bauthor{\bsnm{Liang},~\bfnm{Hua}\binits{H.}}
(\byear{2000}).
\btitle{Quantile regression estimates for a class of linear and
partially linear errors-in-variables models}.
\bjournal{Statist. Sinica}
\bvolume{10}
\bpages{129--140}.
\bid{issn={1017-0405}, mr={1742104}}
\end{barticle}
%

\bptok{imsref}%
% NOT OUTPUTTED:
% number = 1
% fjournal = Statistica Sinica
\endbibitem

%b13 ###
\bibitem{Heiler}
%
\begin{barticle}[mr]
\bauthor{\bsnm{Heiler},~\bfnm{Siegfried}\binits{S.}} \AND
\bauthor{\bsnm{Willers},~\bfnm{Reinhart}\binits{R.}}
(\byear{1988}).
\btitle{Asymptotic normality of {R}-estimates in the linear model}.
\bjournal{Statistics}
\bvolume{19}
\bpages{173--184}.
\bid{doi={10.1080/02331888808802084}, issn={0233-1888}, mr={0945375}}
\end{barticle}
%

\bptok{imsref}%
% NOT OUTPUTTED:
% number = 2
% doi = http://dx.doi.org/10.1080/02331888808802084
% fjournal = Statistics. A Journal of Theoretical and Applied Statistics
\endbibitem

%b14 ###
\bibitem{Hodges}
%
\begin{barticle}[mr]
\bauthor{\bsnm{Hodges},~\bfnm{J.~L.}\binits{J.L.} \bsuffix{Jr.}} \AND
\bauthor{\bsnm{Lehmann},~\bfnm{E.~L.}\binits{E.L.}}
(\byear{1963}).
\btitle{Estimates of location based on rank tests}.
\bjournal{Ann. Math. Statist.}
\bvolume{34}
\bpages{598--611}.
\bid{issn={0003-4851}, mr={0152070}}
\end{barticle}
%

\bptok{imsref}%
% NOT OUTPUTTED:
% fjournal = Annals of Mathematical Statistics
\endbibitem

%b15 ###
\bibitem{Hyk}
%
\begin{barticle}[pbm]
\bauthor{\bsnm{Hyk},~\bfnm{Wojciech}\binits{W.}} \AND
\bauthor{\bsnm{Stojek},~\bfnm{Zbigniew}\binits{Z.}}
(\byear{2013}).
\btitle{Quantifying uncertainty of determination by standard additions
and serial dilutions methods taking into account standard uncertainties
in both axes}.
\bjournal{Anal. Chem.}
\bvolume{85}
\bpages{5933--5939}.
\bid{doi={10.1021/ac4007057}, issn={1520-6882}, pmid={23678943}}
\end{barticle}
%

\bptok{imsref}%
% NOT OUTPUTTED:
% number = 12
% fjournal = Analytical chemistry
\endbibitem

%b16 ###
\bibitem{Hyslop}
%
\begin{barticle}[mr]
\bauthor{\bsnm{Hyslop},~\bfnm{Dean~R.}\binits{D.R.}} \AND
\bauthor{\bsnm{Imbens},~\bfnm{Guido~W.}\binits{G.W.}}
(\byear{2001}).
\btitle{Bias from classical and other forms of measurement error}.
\bjournal{J.~Bus. Econom. Statist.}
\bvolume{19}
\bpages{475--481}.
\bid{doi={10.1198/07350010152596727}, issn={0735-0015}, mr={1963378}}
\end{barticle}
%

\bptok{imsref}%
% NOT OUTPUTTED:
% number = 4
% doi = http://dx.doi.org/10.1198/07350010152596727
% fjournal = Journal of Business \& Economic Statistics
\endbibitem

%b17 ###
\bibitem{Jaeckel}
%
\begin{barticle}[mr]
\bauthor{\bsnm{Jaeckel},~\bfnm{Louis~A.}\binits{L.A.}}
(\byear{1972}).
\btitle{Estimating regression coefficients by minimizing the dispersion
of the residuals}.
\bjournal{Ann. Math. Statist.}
\bvolume{43}
\bpages{1449--1458}.
\bid{issn={0003-4851}, mr={0348930}}
\end{barticle}
%

\bptok{imsref}%
% NOT OUTPUTTED:
% fjournal = Annals of Mathematical Statistics
\endbibitem

%b18 ###
\bibitem{Jur1969}
%
\begin{barticle}[mr]
\bauthor{\bsnm{Jure{\v{c}}kov{\'a}},~\bfnm{Jana}\binits{J.}}
(\byear{1969}).
\btitle{Asymptotic linearity of a rank statistic in regression parameter}.
\bjournal{Ann. Math. Statist.}
\bvolume{40}
\bpages{1889--1900}.
\bid{issn={0003-4851}, mr={0248931}}
\end{barticle}
%

\bptok{imsref}%
% NOT OUTPUTTED:
% fjournal = Annals of Mathematical Statistics
\endbibitem

%b19 ###
\bibitem{Jur1971}
%
\begin{barticle}[mr]
\bauthor{\bsnm{Jure{\v{c}}kov{\'a}},~\bfnm{Jana}\binits{J.}}
(\byear{1971}).
\btitle{Nonparametric estimate of regression coefficients}.
\bjournal{Ann. Math. Statist.}
\bvolume{42}
\bpages{1328--1338}.
\bid{issn={0003-4851}, mr={0295487}}
\end{barticle}
%

\bptok{imsref}%
% NOT OUTPUTTED:
% fjournal = Annals of Mathematical Statistics
\endbibitem

%b20 ###
\bibitem{Jurecko1}
%
\begin{barticle}[mr]
\bauthor{\bsnm{Jure{\v{c}}kov{\'a}},~\bfnm{Jana}\binits{J.}},
\bauthor{\bsnm{Picek},~\bfnm{Jan}\binits{J.}} \AND
\bauthor{\bsnm{Saleh},~\bfnm{A.~K.~Md.~Ehsanes}\binits{A.K.Md.E.}}
(\byear{2010}).
\btitle{Rank tests and regression and rank score tests in measurement
error models}.
\bjournal{Comput. Statist. Data Anal.}
\bvolume{54}
\bpages{3108--3120}.
\bid{doi={10.1016/j.csda.2009.08.020}, issn={0167-9473}, mr={2727738}}
\end{barticle}
%

\bptok{imsref}%
% NOT OUTPUTTED:
% number = 12
% doi = http://dx.doi.org/10.1016/j.csda.2009.08.020
% fjournal = Computational Statistics \& Data Analysis
\endbibitem

%b21 ###
\bibitem{Kelly}
%
\begin{barticle}[auto:parserefs-M02]
\bauthor{\bsnm{Kelly},~\bfnm{B.~C.}\binits{B.C.}}
(\byear{2007}).
\btitle{Some aspects of measurement error in linear regression of
astronomical data}.
\bjournal{The Astrophysical Journal}
\bvolume{665}
\bpages{1489--1506}.
\end{barticle}
%

\bptok{imsref}%
\endbibitem

%b22 ###
\bibitem{Koul2000}
%
\begin{bbook}[mr]
\bauthor{\bsnm{Koul},~\bfnm{Hira~L.}\binits{H.L.}}
(\byear{2002}).
\btitle{Weighted Empirical Processes in Dynamic Nonlinear Models}.
\bseries{Lecture Notes in Statistics}
\bvolume{166}.
\blocation{New York}:
\bpublisher{Springer}.
%\bnote{Second edition of {{\i}t Weighted empiricals and linear models}
%[Inst. Math. Statist., Hayward, CA, 1992; MR1218395 (95c:62061)]}.
\bid{doi={10.1007/978-1-4613-0055-7}, mr={1911855}}
\end{bbook}
%

\bptok{imsref}%
% NOT OUTPUTTED:
% doi = http://dx.doi.org/10.1007/978-1-4613-0055-7
% isbn = 0-387-95476-7
% fpage = xviii+425
\endbibitem

%b23 ###
\bibitem{Marques}
%
\begin{barticle}[mr]
\bauthor{\bsnm{Marques},~\bfnm{Tiago~A.}\binits{T.A.}}
(\byear{2004}).
\btitle{Predicting and correcting bias caused by measurement error in
line transect sampling using multiplicative error models}.
\bjournal{Biometrics}
\bvolume{60}
\bpages{757--763}.
\bid{doi={10.1111/j.0006-341X.2004.00226.x}, issn={0006-341X}, mr={2089452}}
\end{barticle}
%

\bptok{imsref}%
% NOT OUTPUTTED:
% number = 3
% doi = http://dx.doi.org/10.1111/j.0006-341X.2004.00226.x
% fjournal = Biometrics. Journal of the International Biometric Society
\endbibitem

%b24 ###
\bibitem{Ivo}
%
\begin{bmisc}[auto:parserefs-M02]
\bauthor{\bsnm{M{\"{u}}ller},~\bfnm{I.}\binits{I.}}
(\byear{1996}).
\bhowpublished{Robust methods in the linear calibration model.
Ph.D. thesis, Charles Univ. in Prague}.
\end{bmisc}
%

\bptok{imsref}%
\endbibitem

%b25 ###
\bibitem{Navratil}
%
\begin{binproceedings}[auto:parserefs-M02]
\bauthor{\bsnm{Navr{\'{a}}til},~\bfnm{R.}\binits{R.}}
(\byear{2012}).
\btitle{Rank Tests and R-estimates in Location Model with Measurement errors}.
In \bbooktitle{Proceedings of Workshop of the Jaroslav H\'{a}jek Center
and Financial Mathematics in Practice I. Book of Short Papers}
(\beditor{\bfnm{J.}\binits{J.}~\bsnm{Zelinka}}
\AND
\beditor{\bfnm{J.}\binits{J.}~\bsnm{Horov\'a}}, eds.).
\blocation{Brno}:
\bpublisher{Masaryk Univ.}
\end{binproceedings}
%

\bptok{imsref}%
\endbibitem

%b26 ###
\bibitem{NavratilSaleh}
%
\begin{barticle}[mr]
\bauthor{\bsnm{Navr{\'a}til},~\bfnm{Radim}\binits{R.}}
\AND
\bauthor{\bsnm{Saleh},~\bfnm{A.~K.~Md.~Ehsanes}\binits{A.K.Md.E.}}
(\byear{2011}).
\btitle{Rank tests of symmetry and R-estimation of location parameter
under measurement errors}.
\bjournal{Acta Univ. Palack. Olomuc. Fac. Rerum Natur. Math.}
\bvolume{50}
\bpages{95--102}.
\bid{issn={0231-9721}, mr={2920711}}
\end{barticle}
%

\bptok{imsref}%
% NOT OUTPUTTED:
% number = 2
% fjournal = Acta Universitatis Palackianae Olomucensis. Facultas Rerum
%Naturalium. Mathematica
\endbibitem

%b27 ###
\bibitem{Oosterhoff79}
%
\begin{bincollection}[mr]
\bauthor{\bsnm{Oosterhoff},~\bfnm{J.}\binits{J.}} \AND
\bauthor{\bparticle{van} \bsnm{Zwet},~\bfnm{W.~R.}\binits{W.R.}}
(\byear{1979}).
\btitle{A note on contiguity and {H}ellinger distance}.
In \bbooktitle{Contributions to Statistics. Jaroslav H\'ajek Memorial Volume}
(\beditor{\bfnm{J.}\binits{J.}~\bsnm{Jure\v{c}kov\'a}}, ed.)
\bpages{157--166}.
\blocation{Dordrecht}:
\bpublisher{Reidel}.
\bid{mr={0561267}}
\end{bincollection}
%

\bptok{imsref}%
\endbibitem

%b28 ###
\bibitem{Picek}
%
\begin{bmisc}[auto:parserefs-M02]
\bauthor{\bsnm{Picek},~\bfnm{J.}\binits{J.}}
(\byear{1996}).
\bhowpublished{Statistical procedures based on regression rank scores.
Ph.D. thesis, Charles Univ. in Prague}.
\end{bmisc}
%

\bptok{imsref}%
\endbibitem

%b29 ###
\bibitem{Pollard91}
%
\begin{barticle}[mr]
\bauthor{\bsnm{Pollard},~\bfnm{David}\binits{D.}}
(\byear{1991}).
\btitle{Asymptotics for least absolute deviation regression estimators}.
\bjournal{Econometric Theory}
\bvolume{7}
\bpages{186--199}.
\bid{doi={10.1017/S0266466600004394}, issn={0266-4666}, mr={1128411}}
\end{barticle}
%

\bptok{imsref}%
% NOT OUTPUTTED:
% number = 2
% doi = http://dx.doi.org/10.1017/S0266466600004394
% fjournal = Econometric Theory
\endbibitem

%b30 ###
\bibitem{Rocke}
%
\begin{barticle}[auto:parserefs-M02]
\bauthor{\bsnm{Rocke},~\bfnm{D.~M.}\binits{D.M.}} \AND
\bauthor{\bsnm{Lorenzato},~\bfnm{S.}\binits{S.}}
(\byear{1995}).
\btitle{A two-component model for measurement error in analytical chemistry}.
\bjournal{Technometrics}
\bvolume{37}
\bpages{176--184}.
\end{barticle}
%

\bptok{imsref}%
% NOT OUTPUTTED:
% number = 2
\endbibitem

%b31 ###
\bibitem{Saleh2}
%
\begin{barticle}[mr]
\bauthor{\bsnm{Saleh},~\bfnm{A.~K.~Md.~Ehsanes}\binits{A.K.Md.E.}},
\bauthor{\bsnm{Picek},~\bfnm{Jan}\binits{J.}} \AND
\bauthor{\bsnm{Kalina},~\bfnm{Jan}\binits{J.}}
(\byear{2012}).
\btitle{R-estimation of the parameters of a multiple regression model
with measurement errors}.
\bjournal{Metrika}
\bvolume{75}
\bpages{311--328}.
\bid{doi={10.1007/s00184-010-0328-2}, issn={0026-1335}, mr={2909549}}
\end{barticle}
%

\bptok{imsref}%
% NOT OUTPUTTED:
% number = 3
% doi = http://dx.doi.org/10.1007/s00184-010-0328-2
% coden = MTRKA8
% fjournal = Metrika. International Journal for Theoretical and Applied
%Statistics
\endbibitem

%b32 ###
\bibitem{JISA}
%
\begin{barticle}[mr]
\bauthor{\bsnm{Sen},~\bfnm{Pranab~Kumar}\binits{P.K.}},
\bauthor{\bsnm{Jure{\v{c}}kov{\'a}},~\bfnm{Jana}\binits{J.}} \AND
\bauthor{\bsnm{Picek},~\bfnm{Jan}\binits{J.}}
(\byear{2013}).
\btitle{Rank tests for corrupted linear models}.
\bjournal{J. Indian Statist. Assoc.}
\bvolume{51}
\bpages{201--229}.
\bid{issn={0537-2585}, mr={3234614}}
\end{barticle}
%

\bptok{imsref}%
% NOT OUTPUTTED:
% number = 1
% fjournal = Journal of the Indian Statistical Association
\endbibitem

\end{thebibliography}
\end{document}